\documentclass[11pt]{article}

\usepackage{amssymb,amsmath, pinlabel, mathabx, amscd, tikz, amsthm, xcolor,amsfonts}
\usepackage[normalem]{ulem} 
\usepackage[pagebackref]{hyperref}
\hypersetup{colorlinks=true, linkcolor=cyan,citecolor=cyan, urlcolor= cyan}
\usepackage[shortlabels]{enumitem}
\usepackage{mathtools}
\usepackage{lineno}
\input xy
\xyoption{all}

\usepackage[latin1]{inputenc}
\usepackage[T1]{fontenc}
\usepackage{amsfonts,amssymb}
\usepackage[all]{xy}
\usepackage{comment} 
\usepackage{enumitem}
\usepackage{marginnote}

\definecolor{green}{rgb}{0.0, 0.5, 0.0}
\definecolor{yellow}{rgb}{0.8, 0.33, 0.0}

\setlist[itemize,enumerate]{leftmargin=*}

\newcommand{\bb}{\bar b}
\newcommand{\Z}{\mathbb{Z}}

\newcommand{\N}{\mathbb{N}}

\newcommand{\ord}{\mathop{\mathrm{ord}}}
\newcommand{\Zer}{\mathrm{Zer}}
\newcommand{\set}[1]{\left\{#1\right\}}

\newcommand{\de}{\mathbf{i}}
\newcommand{\ex}{\mathbf{e}}
\newcommand{\ic}{\mathbf{c}}

\numberwithin{equation}{section}

\newtheorem{theorem}{Theorem}[section]
\newtheorem{lemma}[theorem]{Lemma}
\newtheorem{pro}[theorem]{Proposition}
\newtheorem{corollary}[theorem]{Corollary}
\newtheorem{property}[theorem]{Property}
\newtheorem{properties}[theorem]{Properties}
\theoremstyle{definition}
\newtheorem{definition}[theorem]{Definition}

\newtheorem{remark}[theorem]{Remark}
\newtheorem{example}[theorem]{Example}

\usepackage{todonotes}

\title{Eggers decomposition of polar curves in positive characteristic}
\author{Ana Bel\'en de Felipe, Evelia R. Garc\'{i}a Barroso, Janusz Gwo\'zdziewicz}
\date{}

\begin{document}


\maketitle

\begin{abstract}
    Let $f$ be a power series in two variables with coefficients in an algebraically closed field of positive characteristic. We identify a necessary and sufficient condition, read from the Eggers-Wall tree of $f$, for $\frac{\partial f}{\partial y}$ to admit Eggers decomposition.
\end{abstract}

\bigskip


\medskip


\section{Introduction}

Let us consider a reduced germ of a plane curve $C$ defined by the equation 
$f(x, y) = 0$ in $(\mathbb{C}^2, 0)$. Let \( \ell(x, y)=0 \) be the equation of a line passing through the origin.
Given this setup, we define the morphism $\Phi: (\mathbb{C}^2, 0) \rightarrow (\mathbb{C}^2, 0)$ by $\Phi(x, y) = (u, v) = (\ell(x, y),f(x, y)).$
The critical locus $P_{\ell}(f)$ of $\Phi$ is called the {\it polar curve} of $C$ with respect to $\ell$. The curve $P_{\ell}(f)$  is the set of points where the tangent to the level curve $f(x, y) =\lambda$ ($\lambda \in \mathbb C$) is parallel to the line given by $\ell$, together with the origin.  According to Teissier (\cite{Teissier1977}), there exists a Zariski open set of the complex projective line for which the polar curve $P_{\ell}(f)$  is reduced and has multiplicity one less than the multiplicity of the germ $f(x,y)=0$.  A polar curve $P_{\ell}(f)$ with  $\ell$ contained in this Zariski open set is referred to as a {\it generic} polar curve.
Observe that  $P_{\ell}(f)$ admits as an equation the jacobian of $\Phi$, and in particular if $\ell=x$ then $P_{\ell}(f)$ is given by the  partial derivative $\frac{\partial f}{\partial y}$ of $f$ with respect to $y$.

The {\it discriminant} of $\Phi$ is the image of $P_{\ell}(f)$ by $\Phi$ and it lives in the complex plane of coordinates $u,v$. The Newton polygon of this discriminant in the coordinates
$u,v$, for a general choice of $\ell$, is independent of $\ell$ and is called the {\it Jacobian Newton polygon} of $f$. This polygon was introduced by  Teissier in \cite{Teissier1977} for isolated hypersurface singularities. Teissier shows that, under these more general settings, the Jacobian Newton polygon is an invariant of the equisingularity, although he does not determine it explicitly. In \cite{Merle1977}, Merle computes it explicitly  for an irreducible power series $f(x, y)$ and proves that it is a complete invariant for the equisingularity class of the germ $C$ defined by $f(x, y) = 0$; that is,  it can be determined from the equisingularity class of  $C$ (from its Puiseux exponents or, equivalently,  from its semigroup of values)  and, in turn, fully characterizes it.

Merle's result answers an original question posed by Brieskorn (see \cite[page 401]{Brieskorn}): how to effectively obtain the equisingularity class of an irreducible plane curve germ from its Jacobian ideal? Merle's strategy is to study the {\it contact} of $C$ with the irreducible components of a generic polar $P_{\ell}(f)$, providing a decomposition of the power series defining $P_{\ell}(f)$ as a product of power series that are not necessarily irreducible, but where the decomposition (number of factors, orders of each factor and their intersection multiplicities with $f$) depends exclusively on the equisingularity class of $C$.  This is very interesting because, as Pham illustrates with his example (see \cite[Exemple 3]{Pham}),  the equisingularity class of the generic polar curves $P_{\ell}(f)$ depends not only on the equisingularity class of the germ of plane curve $C$ but also on its analytic type. However, Merle's decomposition captures the information of $P_{\ell}(f)$ that depends only on the equisingularity of $C$. 

Eggers, in \cite{Eggers}, presents an example illustrating that the Jacobian Newton polygon of $f$ fails to be a complete invariant when 
$f$ is a reduced, non-irreducible power series. Eggers represents the equisingularity class of $f$ by means of a tree, that the second-named author revisits in \cite{GBtesis} and names the Eggers tree, 
and using it determines a decomposition of the polar curve $P_{\ell}(f)$ in the sense of Merle. We will henceforth refer to this decomposition as the {\it Eggers decomposition} of the polar curve, to emphasize that it is read from the tree. In \cite{Walltree}, Wall modified these Eggers trees slightly, and they became known in the literature as Eggers-Wall trees.

In \cite{GB-P}, the authors extend Merle's decomposition of $\frac{\partial f}{\partial y}$ to the case where $f$ is an irreducible power series in two variables over an algebraically closed field $K$ of positive characteristic. They obtain a decomposition similar to that given by Merle for $K=\mathbb C$, but assuming additionally
that the characteristic of $K$ does not divide the order of $f(0,y)$. The initial purpose of this work was to investigate under what hypothesis we obtain an Eggers decomposition of  $\frac{\partial f}{\partial y}$ for the case when $f$ is  non-irreducible and $K$ is algebraically closed of positive characteristic, that is, under what assumptions we can generalise to positive characteristic the decomposition of the polar curve given in \cite{GBtesis} (see also \cite{GB 2000}).

The reduced non-irreducible case in arbitrary characteristic is considerably more delicate, as we show in this paper. We explicitly detect a condition that the Eggers-Wall tree of $f$ must verify in order to guarantee an Eggers decomposition of  $\frac{\partial f}{\partial y}$. We call this condition the {\it Eggers condition}. We get this Eggers decomposition  (Corollary \ref{Eggersfact}) as a consequence of the main theorem (Theorem \ref{mainJ}) of the paper. Moreover we prove that Eggers condition is a necessary and sufficient condition for obtaining Eggers decomposition of  $\frac{\partial f}{\partial y}$. Meeting Eggers condition requires first passing through other conditions ((E1), (E2) and $\de$ conditions) that are also interpreted in the Eggers-Wall tree of $f$.

Before reaching these conditions, we must overcome certain issues. The first one arises in determining which Eggers-Wall tree we are going to consider in arbitrary characteristic (see \cite[Section 16]{Jaca}). Eggers originally defined his tree using Puiseux parametrisations of the curve germ $C$ given by $f(x,y)=0$, but it is well known that these parametrisations do not necessarily exist in positive characteristic, and even if $C$ admits Puiseux parametrisation, we do not necessarily have as many as in the case of characteristic zero. On the other hand, as Example \ref{ex:pdoesnotimplyE} illustrates, the Eggers condition on $f$ does not force the existence of a Puiseux factorization of $\frac{\partial f}{\partial y}$, so in order to prove the Eggers decomposition of $\frac{\partial f}{\partial y}$ we cannot apply the classical Puiseux roots approach, as in characteristic zero, and we need to use another set of tools.

From now on we denote by $K$ an algebraically closed field of characteristic $p\geq0$. The paper is organized as follows. In Section \ref{sect2}, we review fundamental concepts related to plane curve singularities in arbitrary characteristic. Specifically, we recall the notions of logarithmic distance, the order of coincidence,  the Eggers-Wall tree and its associated functions, the semigroup of values and the Puiseux characteristic sequence (in the irreducible case) and the Newton diagram. We also introduce what we call the Eggers condition along with the $\de$-condition, which will be essential in the development of this work. In Section \ref{sec:Ed}, for a given pair of power series $(f,\Gamma)$, we introduce the conditions (E1) and (E2), which extract information from the second power series $\Gamma$ by means of the Eggers-Wall tree $\Theta(f)$ of $f$. We then prove (Theorem \ref{thm:E1iffE2}) that these two conditions are equivalent.  Condition (E1) involves expressing  $\Gamma$ as a  product of power series, not necessarily irreducible, and the value of the intersection multiplicity of each factor of this product with the power series $x$. For this reason, under the hypothesis $0<\ord f(0,y)<\infty$, we say that $\Gamma$ admits an Eggers decomposition with respect to $f$ if either of these equivalent conditions (E1) or (E2) is satisfied by the pair  $(\Gamma,f)$. In particular, if $1<\ord f(0,y)<\infty$ then we say that $\frac{\partial f}{\partial y}$ \textit{admits Eggers decomposition} if $\frac{\partial f}{\partial y}$ admits Eggers decomposition with respect to $f$. We conclude Section \ref{sec:Ed} by providing information on the intersection multiplicity of $x$ and the irreducible factors of $f$
with the factors of the Eggers decomposition of $\Gamma$ (Corollary \ref{coro:Edecg}). This information only depends on the  Eggers-Wall tree $\Theta(f)$. Moreover, we compute the exponent function associated with the  Eggers-Wall tree  $\Theta(f)$,  using the factors of the Eggers decomposition of $\Gamma$ and the irreducible factors of $f$ (Theorem \ref{thm:1}). In Section \ref{sec:sufficient}, after proving two technical lemmas, we establish the main theorem of this paper (Theorem \ref{mainJ}), and as a consequence, we conclude  that if $f$ satisfies Eggers condition, then $\frac{\partial f}{\partial y}$ admits Eggers decomposition (Corollary \ref{Eggersfact}). We end this fourth section with several illustrative examples. Finally, in Section \ref{sec:5}, under the assumption that $\frac{\partial f}{\partial y}$
admits Eggers decomposition, we prove that $f$ satisfies Eggers condition (Theorem \ref{necessary}). This is achieved by first establishing that $\de$-condition holds for $f$ as an intermediate step (Proposition \ref{i-condition}). At this point, we would like to emphasize that the statement of Proposition \ref{i-condition} is not an equivalence as Example \ref{ex:pdoesnotimplyE} illustrates. The proof of Proposition \ref{i-condition} relies  on the fact that the shape of the Eggers-Wall tree of 
$f$ is preserved under the ramification induced by replacing in $f$ the variable $x$
with some power of $x$ whose exponent involves the characteristic $p$ of the field 
$K$ (Remark \ref{rem:hat}).  A series of lemmas at the beginning of Section \ref{sec:5} are devoted to showing this.


\section{Background on plane curve singularities}
\label{sect2}

For any power series $f,g\in K[[x,y]]$, 
their intersection multiplicity at the origin is $i_0(f,g):=\dim_K K[[x,y]]/(f,g)$,
where $(f,g)$ is the ideal of $K[[x,y]]$ generated by $f$ and $g$. If $f$ and $g$ are both irreducible and coprime with $x$, we put
\[
d(f,g):=\frac{i_0(f,g)}{i_0(f,x)i_0(g,x)}.
\]
The function $d$ is a logarithmic distance (see \cite[Theorem 2.8]{GB-P-2015}). In particular, 
$d$ satisfies the so-called {\it strong triangle inequality} (STI): for any three irreducible power series $f,g,$ and $h$ coprime with $x$,
\[
d(f,g)\geq \inf\set{d(f,h), d(g,h)},
\]
which is equivalent that at least two of the numbers $d(f,g), d(f,h), d(g,h)$ are equal and the third is no smaller than the other two.

A {\it Puiseux series} is a power series in the ring $\bigcup_{n\in\N\setminus\set{0}} K[[x^{1/n}]]$. 
The \textit{support} of a Puiseux series $\alpha(x)=\sum_{i\geq0} a_ix^{i/n}$ is the set of rational numbers $i/n$ such that $a_i\neq0$. 
The {\it index} of $\alpha(x)$ is the minimal common denominator of the elements of its support.

\begin{property}\label{prop:minpol}
    The minimal polynomial of a Puiseux series $\alpha(x)=\sum_{i\geq0} a_ix^{i/n}$ over the field of fractions of $K[[x]]$ is an irreducible Weierstrass polynomial of degree equal to the index of $\alpha(x)$.
\end{property}

\begin{proof}
By \cite[Theorem 2.6]{Pl2013}, there exists an irreducible power series $f\in K[[x,y]]$ such that $f(t^n,\alpha(t^n))=0$. 
By the Weierstrass  Preparation Theorem, $f$ is a product of a unit of the ring $K[[x,y]]$ and an irreducible Weierstrass 
polynomial $w\in K[[x]][y]$. We have that $w(x,\alpha(x))=0$. Since $K[[x]]$ is a unique factorization domain, 
Gauss theorem tell us that $w$, when treated as a polynomial with coefficients in the field of fractions of $K[[x]]$, remains irreducible. 
Therefore, $w$ is the minimal polynomial of $\alpha(x)$ over the fraction field of $K[[x]]$. The statement about the degree follows from \cite[Proposition 17.5--(2)]{Jaca}.
\end{proof}

Let $f\in K[[x]][y]$ be an irreducible Weierstrass polynomial of degree $d$. A {\it Puiseux root} of $f$ is any Puiseux series $\alpha(x)$ such that $f(x,\alpha(x))=0$. Denote by $\Zer(f)$ the set of \textit{Puiseux roots} of $f$ (it could be the empty set). If $p=0$ or if $p$ does not divide $d$, then $f$ has $d$ distinct roots in $K[[x^{1/d}]]$ (this is the Newton--Puiseux theorem in characteristic zero, and in positive characteristic  see \cite[Proposition 2.1.8, Chapter II]{Campillo} and  \cite[Corollary 2.4]{Pl2013}). The {\it order of coincidence} of $f$ and $g$, where $g$ is also an irreducible Weierstrass polynomial in $K[[x]][y]$ and the sets of Puiseux roots of $f$ and $g$ are non-empty, is defined as
\[k(f,g):=\max\set{\ord(\alpha(x)-\beta(x))\;:\;\alpha(x)\in\Zer(f) \;\text{and}\; \beta(x)\in\Zer(g) }.\]

For any irreducible power series $f\in K[[x,y]]$, the {\it semigroup} associated with the branch defined by $f$ is $\Gamma(f):=\set{i_0(f,g)\;:\;g\notequiv 0 \pmod f}$. 

\subsection{The Eggers-Wall tree}
\label{sect:E-Wtree}
In this subsection, we consider a power series $f$ in $K[[x,y]]$ with $f(0,0)=0$. We introduce the {\it Eggers-Wall tree} associated with $f$ with respect to $x$. Since the power series $x$ will be fixed throughout this article, we will call it the Eggers-Wall tree of $f$. We also present related notions that are used elsewhere in the paper. In particular, we introduce the {\it Eggers condition} and the $\de$--{\it condition} which are essential in this work. For further details about trees, we refer to \cite[Section 1.6]{Handbook}, where the Eggers-Wall tree of $f$ is introduced using $\Zer(f)$ instead of the semigroups of irreducible components of $f$, but in characteristic zero both trees are equivalent. Here, we use the semigroups instead of $\Zer(f)$,  because the latter set may be empty, while the semigroup of an irreducible series behaves similarly in any characteristic (see \cite{Jaca} for details on the Eggers-Wall tree in arbitrary characteristic).

In what follows, we assume that $f$ is irreducible and coprime with $x$. Let $\bb_0,\ldots,\bb_h\in\Gamma(f)$ be such that $\bb_0=i_0(f,x)$, $\bb_i=\min(\Gamma(f)\setminus(\N\bb_0+\ldots+\N\bb_{i-1}))$ for $i\in\set{1,\ldots,h}$, and $\Gamma(f)=\N\bb_0+\ldots+\N\bb_h$ (such a sequence exists and  depends only on $\bb_0$ and $\Gamma(f)$, and not on the characteristic of $K$, see \cite[Theorem 4.3.8, Chapter IV]{Campillo} or \cite[Chapter 6, Proposition 6.1]{He2003} ). Set $e_i=\gcd(\bb_0,\ldots,\bb_i)$ for $i\in\set{0,\ldots,h}$. Then $e_0>e_1>\ldots>e_h=1$ and $n_i \bb_i<\bb_{i+1}$ for $i\in\set{1,\ldots,h-1}$, where $n_i:=e_{i-1}/e_i$ \label{enei} for $i\in\set{1,\ldots,h}$. Moreover, there exists a sequence of monic irreducible polynomials $f^{(0)},\ldots,f^{(h-1)}\in K[[x]][y]$ such that $\deg_y f^{(i)}=i_0(f^{(i)},x)=\bb_0/e_{i}$ and $i_0(f,f^{(i)})=\bb_{i+1}$ for $i\in\set{0,\ldots,h-1}$. We refer the reader to \cite[Section 4]{GB-P-2015} for more details about these results. 

The sequence $f^{(0)},\ldots,f^{(h-1)},f^{(h)}\in K[[x]][y]$, where $f^{(h)}\in K[[x]][y]$ denotes the Weierstrass polynomial associated to $f$, is called a \textit{sequence of key polynomials} of $f$. The polynomial $f^{(i)}$ is called an $i$-th key polynomial of $f$, for $0\leq i\leq h$. A direct computation shows that $d(f^{(0)},f)<d(f^{(1)},f)<\ldots<d(f^{(h-1)},f)<d(f^{(h)},f)=\infty$.  

The Eggers-Wall tree $\Theta(f)$ of $f$ is a compact segment endowed with the following:
\begin{itemize}
    \item An increasing homeomorphism $\ic\colon\Theta(f)\to[0,\infty]$,
    called the \textit{contact complexity function}.
    \item A set of \textit{marked points}, which are by definition the preimages by the function $\ic$ of the values $d(f^{(0)},f),\ldots,d(f^{(h-1)},f)$, as well as $\ic^{-1}(0)$, labelled by $x$, and $\ic^{-1}(\infty)$, labeled by $f$.
    \item A piecewise constant function $\de\colon\Theta(f)\to\N$, called the \textit{index function}, which takes at a point $P\in\Theta(f)$ the value $\bb_0/e_i$, where $i=\min\set{j\;:\;\ic(P)\leq d(f^{(j)},f)}$.
\end{itemize}

\begin{figure}
    \begin{center}
\begin{tikzpicture}[scale=1]
   
    \draw[->, color=black,  line width=1pt](0,1) -- (0,7);
    \draw [->, color=black,  line width=1pt](0,1) -- (0,0);
   
    \node [below, color=black] at (0,0) {$x$};
    \node [right, color=black] at (0,0.5) {\small{$1$}};
    \node [right, color=black] at (0,1.5) {\small{$n_{1}$}};
    \node [right, color=black] at (0,2.5) {\small{$n_{1}n_{2}$}};
    \node [right, color=black] at (0,5.5) {\small{$n_{1}n_{2}\cdots n_{h-1}$}};
    \node [right, color=black] at (0,6.5) {\small{$n_{1}n_{2}\cdots n_{h}$}};
    \node [above, color=black] at (0,7) {$f$};
                        
    \node[draw,circle, inner sep=1.5pt,color=black, fill=black] at (0,1){};
    \node[draw,circle, inner sep=1.5pt,color=black, fill=black] at (0,2){};
    \node[draw,circle, inner sep=1.5pt,color=black, fill=black] at (0,3){};
    \node[draw,circle, inner sep=1.5pt,color=black, fill=black] at (0,4){};
    \node[draw,circle, inner sep=1.5pt,color=black, fill=black] at (0,5){};
    \node[draw,circle, inner sep=1.5pt,color=black, fill=black] at (0,6){};
          
    \node [left, color=black] at (0,1) {$P_{1}$}; 
    \node [left, color=black] at (0,2) {$P_{2}$}; 
    \node [left, color=black] at (0,3) {$P_{3}$};
    \node [left, color=black] at (-0.2,4) {$\vdots$};
    \node [left, color=black] at (0,6) {$P_{h}$};  
    \node [left, color=black] at (0,5) {$P_{h-1}$};  
    \end{tikzpicture}
\end{center}
 \caption{Eggers-Wall tree of an irreducible power series $f$}
\label{fig:EW irred}
\end{figure}
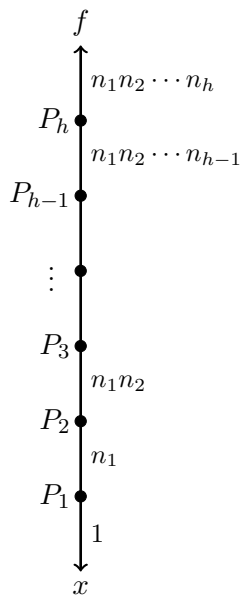

In Figure \ref{fig:EW irred} we have drawn the Eggers-Wall tree $\Theta(f)$ of $f$ , together with the different values attained by the index function. 
For $1\leq i\leq h$, the point $P_i$ is the marked point of $\Theta(f)$ with $\ic(P_i)=d(f^{(i-1)},f)=e_{i-1}\bb_i/\bb_0^2$. 
Observe that if we define $n_0:=1$, then $\bb_0/e_i=n_0\cdots n_i$ for $0\leq i\leq h$. 
Since $e_i=\gcd(\bb_i,e_{i-1})$, the integers $m_i:=\bb_i/e_i$ and $n_i=e_{i-1}/e_i$ are coprime. Then 
\begin{equation}\label{contact-to-index}
    \ic(P_i)=\frac{m_i}{(n_0\cdots n_{i-1})^2 n_i}.
\end{equation}

Two marked points $P,Q$ of $\Theta(f)$ are \textit{consecutive} if $\ic(P)<\ic(Q)$ and there is no other marked point in the open segment $(P,Q)$. 

By convention, the Eggers-Wall tree $\Theta(x)$ is reduced to a single point labeled by $x$, with $\ic(x)=0$ and $\de(x)=1$. Let us consider now the case $f=f_1\cdots f_r$, where $r>1$, $f_i\in K[[x,y]]$ is irreducible for $1\leq i\leq r$, and $f_i$ and $f_j$ are coprime for $i\neq j$. Take a pair of irreducible factors $f_i$ and $f_j$ with $i\neq j$ and $f_i$ and $f_j$ both coprime with $x$. Denote by $\ic_i$ (resp. $\ic_{j}$) the contact complexity function on $\Theta(f_i)$ (resp. on $\Theta(f_j)$). Next, we show that the segments $[x,\ic_i^{-1}(d(f_i,f_j))]$ and $[x,\ic_{j}^{-1}(d(f_i,f_j))]$ of the Eggers-Wall trees of $f_i$ and $f_j$ respectively, coincide. Let $\de_i$ (respectively $\de_j$) be the index function on $\Theta(f_i)$ (resp. on $\Theta(f_j)$).

\begin{lemma}\label{lemma:gluing}
     Let $P_0=x$, $P_1, \dots, P_{k+1}$ (respectively $P_0'=x$, $P_1'$, $\dots P_{l+1}'$) be consecutive marked points of $\Theta(f_i)$ (respectively $\Theta(f_j)$) such that $\ic_i(P_k)< d(f_i,f_j)\leq\ic_i(P_{k+1})$ (respectively $\ic_j(P_l')< d(f_i,f_j)\leq\ic_j(P_{l+1}')$). Then
     \begin{itemize}
     \item $k=l$,
     \item $\ic_i(P_s)=\ic_j(P_s')$ for $1\leq s\leq k$,
     \item $\de_i(P_s)=\de_j(P_s')$ for $1\leq s\leq k+1$.
     \end{itemize}
\end{lemma}

\begin{proof}
If $k=l=0$ then there is nothing to prove. Thus without loss of generality we may assume, switching the role of $f_i$ and $f_j$ if necessary, that $k\geq 1$ and $k\geq l$.

Let $(f_i^{(0)},\dots,f_i^{(k-1)})$ be an initial sequence of key polynomials of $f_i$. Then $d(f_i,f_i^{(s-1)})=\ic_i(P_s)$ for $1\leq s\leq k$. By \cite[Theorem 5.2]{GB-P-2015}, $(f_i^{(0)},\dots,f_i^{(k-1)})$ is also an initial sequence of key polynomials of $f_j$. By the STI we have $d(f_j,f_i^{(s-1)})=d(f_i,f_i^{(s-1)})<d(f_j,f_i)$ for $1\leq s\leq k$. Hence $k=l$ and $\ic_i(P_s)=\ic_j(P_s')$ for $1\leq s\leq k$. By \eqref{contact-to-index} the sequence $\de_i(P_s)$ for $1\leq s\leq k+1$ can be computed recursively as follows:
$\de_i(P_1)=1$, $\de_i(P_{s+1})=n_s \de_i(P_s)$ where $n_s$ is the denominator of the irreducible fraction $\frac{m_s}{n_s}=\ic_i(P_s)\de_i(P_s)^2$. Since this recursive formula depends only on the contact complexity function the last statement of the lemma also follows. 
\end{proof}

Lemma \ref{lemma:gluing} allows us to do the following construction for the reduced series $f=f_1\cdots f_r$ with several irreducible factors $f_{i}$. The Eggers-Wall tree of $\Theta(f)$ is the tree obtained from the disjoint union $\bigsqcup_{i=1}^r\Theta(f_i)$ by the following identification: for every pair of distinct irreducible factors $f_i$ and $f_j$, identify the points $P\in\Theta(f_i)$ and $P'\in\Theta(f_j)$ if $\ic_{i}(P)=\ic_{j}(P')\leq d(f_i,f_j)$. It is endowed with the contact complexity function $\ic\colon\Theta(f_1\cdots f_r)\to[0,\infty]$ and the index function $\de\colon\Theta(f_1\cdots f_r)\to\N$, obtained by gluing the contact complexity functions and the index functions, respectively, of the Eggers-Wall trees $\Theta(f_i)$ for $1\leq i\leq r$. Its set of marked points consists of its ramification points and the images of the marked points of the $\Theta(f_i)$ by the identification map. Observe that the set of marked points of $\Theta(f)$ is the union of the set of ramification points, the set of points of discontinuity of the index function, and the set of endpoints (i.e., the root labeled by $x$, and the leaves labeled by $f_i$ for $1\leq i\leq r$). A {\it bamboo point} of $\Theta(f)$ is a marked point which is neither an endpoint nor a ramification point.  

When the decomposition of $f$ into irreducible factors is of the form $f=f_1^{\alpha_1}\cdots f_r^{\alpha_r}$, with some $\alpha_i$ greater than one, the Eggers-Wall tree of $f$ is the tree $\Theta(f_1\cdots f_r)$. Note that some authors equip $\Theta(f)$ with some additional labels: the leaf labeled by $f_i$ is also labeled by $\alpha_i$, if $\alpha_i>1$. However in this paper we will not consider such labels.

We consider $\Theta(f)$ as a poset with respect to the partial order $\preceq_x$, where $P\preceq_{x}Q$ if the segment $[x,P]$ is contained in the segment $[x,Q]$, for $P,Q\in\Theta(f)$. In this situation $\ic(P)\leq \ic(Q)$. Let $P$ and $Q$ be two marked points of $\Theta(f)$. We say that $P$ and $Q$ are \textit{consecutive} if $P\prec_x Q$ and there are no marked points in $(P,Q)$. If $P$ and $Q$ are consecutive, we call $Q$ a \textit{direct successor} of $P$.

Given three points $P,Q,T\in\Theta(f)$, the \textit{tripod} generated by them is the union of the segments $[P,Q]$, $[Q,T]$, and $[T,P]$. The intersection point of these segments is denoted by $\langle P,Q,T\rangle$ and is called the \textit{center of the tripod}. 
If the three points lie on a segment, say $T\in[P,Q]$, observe that $\langle P,Q,T\rangle=T$. By construction, the \textit{tripod formula} holds, that is,
\[
i_0(f_i , f_j) = \de (f_i )  \de (f_j) \ic ( \langle x, f_i, f_j \rangle ),
\]
for $i,j\in\set{1,\dots,r}$, such that $i\neq j$ and $f_i$ and $f_j$ are both coprime to $x$.

For any irreducible power series $g\in K[[x,y]]$, the \emph{attaching point} of $g$ to the Eggers-Wall tree $\Theta(f)$ is 
\[P_g  := \max \{ \langle x, g , f_i \rangle\;:\;  i =1, \dots, r \} \in \Theta(f),\] 
where the tripods are viewed in $\Theta(fg)$ and the maximum is taken with respect to ${\preceq_{x}}$. If $g$ divides $f$, then we have that $P_g$ is the  endpoint of $\Theta(f)$ labeled by $g$.

There is a third function defined on the Eggers-Wall tree of $f$, $\ex\colon\Theta(f)\to[0,\infty]$, which is called the \textit{exponent function}. Given $P\in\Theta(f)$, the exponent of $P$ is
\begin{equation}\label{def:ex}
\ex(P)=\int_x^P \de\: d\ic.
\end{equation}
This integral can be computed as follows. Choose an irreducible factor $g$ of $f$ such that $P\preceq_x g$ and take a sequence of key polynomials $g^{(0)},\ldots,g^{(h)}$ of $g$. Denote $c_0=0$ and $c_i=d(g^{(i-1)},g)$ for $1\leq i\leq h+1$. For $i\in\set{0,\ldots,h+1}$, let $P_i$ be the marked point in the segment $[x,g]$ such that $\ic(P_i)=c_i$. Observe that the index function is constant on every segment of the form $(P_{i},P_{i+1}]$. If $P=x$ then $\ex(P)=0$. Otherwise, if $k$ is such that $P\in(P_{k},P_{k+1}]$, then $\ex(P)=\sum_{i=1}^k \de(P_i)(c_i-c_{i-1})+\de(P_{k+1})(\ic(P)-c_k)$.

In the case where $f$ is irreducible and $i_0(f,x)$ is not divisible by $p$, we can consider the \textit{Puiseux characteristic sequence} $b_0,b_1,\ldots,b_h$ of $f$ (see \cite[Theorem 3.5.1, Chapter III]{Campillo} or \cite[Sections 1 and  4 of Chapter 3]{P-P} for details). Recall that one usually denotes $b_{h+1}=\infty$. Taking into account the relation between the $\bb_j$ and the $b_j$, that is, $\bb_i=b_i$ for $i\in \{0,1\}$ and $\bb_i=n_{i-1}\bb_{i-1}+b_i-b_{i-1}$ for $i\in \{2,\ldots, h\}$ (see \cite[Proposition 4.3.5, Chapter IV]{Campillo} or \cite[Section 1 of Chapter 3]{P-P}), a direct computation shows that $\ex(P_{i})=\int_x^{P_{i}} \de\: d\ic=b_{i}/b_0$, for $1\leq i\leq h$.

Note that \begin{equation}\label{eq:gcd}e_l:=\gcd(\bar b_0, \ldots, \bar b_l)=\gcd(b_0, \ldots, b_l)\;\;\; \hbox{\rm for any $l\in \{0,\ldots, h\}$}.
\end{equation}

\begin{property}\label{prop:1}
Assume that $f,g\in K[[x]][y]$ are irreducible Weierstrass polynomials with $\deg f\not\equiv 0\pmod p$ and $\deg g \not\equiv 0\pmod p$. If $P$ is the ramification point of $\Theta(fg)$, then $\ex(P)=k(f,g)$.
\end{property}

\begin{proof}
Let us keep the notation introduced in this section related to the tree $\Theta(f)$ for $f$ irreducible. The assumption on the degree of $f$ allows considering the Puiseux characteristic sequence $b_0,b_1,\ldots,b_h$ of $f$. We have that $\ex(P_{i})=\int_x^{P_{i}} \de\: d\ic=b_{i}/b_0$, for $1\leq i\leq h$. Next observe from the definitions that the functions $\ex$ and $\de$ on $\Theta(f)$ determine the contact complexity function. Indeed, for all $P\in\Theta(f)$, 
\begin{equation}\label{integralform}
    \ic(P)=\int_x^P \frac{d\ex}{\de} . 
\end{equation}

Let $P$ be the attaching point of $g$ to $\Theta(f)$ and let $k\in\set{1,\ldots,h+1}$ be the smallest integer such that $k(f,g)\leq b_{k}/b_0$. According to the Noether formula (see \cite[Theorem 4.1.6]{Wall} for $K=\mathbb C$ but which also holds for any characteristic)  $d(f,g)=\sum_{i=1}^{k-1}\frac{e_{i-1}-e_i}{b_0}\frac{b_i}{b_0}+\frac{e_{k-1}}{b_0}k(f,g)$. 
The right-hand side of this equation can be expressed in integral form, obtaining that $\ic(P)=d(f,g)=\int_x^{\ex^{-1}(k(f,g))}{\frac{d\ex}{\de}}$. Since the contact complexity function is strictly increasing, we deduce that $P=\ex^{-1}(k(f,g))$ by comparing with \eqref{integralform}.
\end{proof}

\begin{definition}
Consider a power series $f\in K[[x,y]]$ with $f(0,0)=0$ and a point $P$ of $\Theta(f)$. Let $h\in K[[x,y]]$ be irreducible. We say that $h$ is {\it above} $P$ in $\Theta(f)$ if $P\preceq_{x} P_h$. 
\end{definition}

Let $h=\prod_{i=1}^mh_i$ be the factorization of $h\in K[[x,y]]$ into irreducible factors. If $P$ is a point of $\Theta(f)$ we put
\begin{equation}
h_P:=\prod h_i,
\end{equation}
where the product runs over $i\in\{1,\ldots,m\}$ such that $h_i$ is above $P$ in $\Theta(f)$ (by convention the empty product is $1$). 

\begin{pro}\label{prop:step}
Let $h\in K[[x,y]]$ with $h(0,0)=0$ and let $P'$ and $P$ be two consecutive marked points of $\Theta(f)$. Then the function $\phi:[P',P]\to\Z$ defined by $\phi(Q)=i_0(h_Q ,x)$, is a weakly decreasing step function which is continuous from the left at every point in the interval $[P',P]$. The points of discontinuity of this function are the points of $[P',P)$ where the irreducible factors of $h$ are attached.
\end{pro}

\begin{proof}
Let $[Q_1,Q_2]$ be any sub-interval of $[P',P]$ that does not contain any attaching points of the irreducible factors of $h$. 
Then for every $Q\in [Q_1,Q_2]$ we have $h_Q=h_{Q_2}$, therefore the function $\phi$ is constant in the interval $[Q_1,Q_2]$. 
This shows that $\phi$ is a step function. We also deduce that $\phi$ is continuous from the left at every point of $[P',P]$.
Furthermore, we see that the discontinuity points of $\phi$ can only occur at points where the factors of $h$ are attached.
Next assume that $Q\in [P',P)$ is the attaching point to $\Theta(f)$ of some irreducible factor $h_i$ of $h$. Consider any $Q_1\in (Q,P]$.
Then $h_i h_{Q_1}$ divides $h_{Q}$ and consequently $i_0(h_Q,x)>i_0(h_{Q_1},x)$. 
This proves that the function $\phi$ is weakly decreasing and $Q$ is a point of discontinuity of $\phi$. 
\end{proof}

\begin{definition}\label{def:conditions}
Let $f\in K[[x,y]]$ be a power series with $f(0,0)=0$. Let $P$ be a marked point of $\Theta(f)$ different from the root. We say that $f$ satisfies \textit{Eggers condition} at the point $P$ if $i_0(f_P,x)\notequiv 0 \pmod p$. We say that $f$ satisfies \textit{$\de$--condition} at the point $P$ if $\de(P)\notequiv 0 \pmod p$.

We say that $f$ satisfies \textit{Eggers condition} (respectively, \textit{$\de$--condition}) if $f$ satisfies Eggers condition (respectively, \textit{$\de$--condition}) at any marked point $P$ of $\Theta(f)$ different from the root. 
\end{definition}

\begin{properties}\label{prop:pE} 
Let $f\in K[[x,y]]$ be a power series with $f(0,0)=0$. Denote by $R$ the root of the Eggers-Wall tree of $f$. Then
\begin{enumerate}
\item[(a)] If $f$ verifies the $\de$--condition at a marked point $P$ of $\Theta(f)$ different from $R$, then $f$ verifies the $\de$--condition at any marked point of the segment $(R,P]$.
\item[(b)] If $f$ verifies Eggers condition at a marked point $P\in \Theta(f)$ different from $R$, then there is an irreducible factor $h$ of $f_P$ such that $f$ verifies the $\de$--condition at any marked point in $(R,P_h]$. In particular, if  $f$ verifies Eggers condition at a marked point $P\in \Theta(f)$ different from $R$, then $f$ verifies the $\de$--condition at $P$.
\item [(c)] If $f$ verifies the Eggers condition, then $f$ verifies the $\de$--condition.
\end{enumerate}
\end{properties}

\begin{proof}
Let $P$ be a marked point of $\Theta(f)$ different from $R$. The statement $(a)$ follows from the fact that $\de(P')$ divides $\de(P)$ for any market point $P'\in (R,P]$.  If $f$ verifies Eggers condition at $P$ then there is at least an irreducible component $h$ of $f_P$ such that $i_0(h,x)=\de(P_h)\notequiv 0 \pmod p$, hence by $(a)$, $f$ verifies the $\de$--condition at any marked point of the segment $(R,P_h]$, and in particular $f$ verifies the $\de$--condition at $P$.
\end{proof}

Statement $(a)$ and $(b)$ of Properties \ref{prop:pE} show that $\de$ and Eggers conditions at a marked point  $P\in \Theta(f)$ have a {\it descending propagating effect} of the function $\de$.

\begin{remark}   
If $f=x^ag$ with $a\geq1$, $g(0,0)=0$, and $x\nmid g$, then $1\leq i_0(f_P,x)<\infty$ for all marked point $P$ of $\Theta(f)$ different from the root. Moreover $f=x$ satisfies the Eggers condition and the $\de$--condition, since there is no marked point in $\Theta(x)$ other than its root.
\end{remark}

If $f\in K[[x,y]]$ is irreducible, then $f$ satisfies the $\de$--condition if and only if $f$ satisfies the Eggers condition. Indeed, after Figure \ref{fig:EW irred} the Eggers condition means that $n_0n_1 \ldots n_h\notequiv 0 \pmod p$ and the $\de$--condition signifies that $n_0\cdots n_i \notequiv 0 \pmod p$ for any $i\in \{0,\ldots, h\}$. Moreover, in the irreducible case the Eggers and $\de$--conditions are equivalent to that $n_i\notequiv 0 \pmod p$, for any $i\in \{0,\ldots, h\}$. Nevertheless  the converse of Property \ref{prop:pE}--(c) does not hold, in general, for non-irreducible series, as the following example shows.

\begin{figure}[h]
\begin{center}
\begin{tikzpicture}[scale=1]
    \draw [-, color=black,  line width=1pt](0,0) -- (0,1);
    \draw [->, color=black,  line width=1pt](0,1) -- (1,2.3);
    \draw [->, color=black,  line width=1pt](0,1) -- (-1,2.3);
    \draw [->, color=black,  line width=1pt](0,1) -- (0,0);
   
    \node [below, color=black] at (0,0) {$x$};
    \node [right, color=black] at (0,0.5) {\small{$1$}};
    \node [above, color=black] at (-1,2.3) {$Q$};
    \node [right, color=black] at (0.45,1.6) {\small{$2$}};
    \node [left, color=black] at (-0.45,1.6) {\small{$1$}};
    \node [above, color=black] at (1,2.3) {$T$};
             
    \node[draw,circle, inner sep=1.5pt,color=black, fill=black] at (0,1){};
          
    \node [left, color=black] at (0,1) {$P$};
\end{tikzpicture}
\end{center}
\caption{Eggers-Wall tree of Example \ref{ex:pdoesnotimplyE}.}
\label{fig:exconditions}
\end{figure}
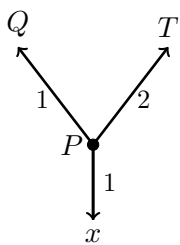

\begin{example}\label{ex:pdoesnotimplyE}
    Consider $f(x,y)=y(y^{2}+y^{3}+x^{5})\in K[[x,y]]$, where $K$ is an algebraically closed field of characteristic $p=3$. The  Eggers-Wall tree of $f$ is drawn on Figure \ref{fig:exconditions}, where $Q$ (respectively $T$) is the leaf corresponding to $y$ (respectively to $y^{2}+y^{3}+x^{5}$). Notice that  $f_{P}=f$, $f_{Q}=y$, $f_{T}=y^{2}+y^{3}+x^{5}$, $\de(P)=\de(Q)=1$, and $\de(T)=2$. At the point $P$, the power series $f$ satisfies the $\de$--condition but not the Eggers condition.  
\end{example}

\subsection{Newton diagrams}

Let $\mathbb R^{+}=\{x \in \mathbb R\;:\; x \geq 0\}$. A subset $N$ of $(\mathbb R^{+})^{2}$ is a {\it Newton diagram} if $N$ is the convex hull of $S+(\mathbb R^{+})^{2}$, for some $S\subset \mathbb N^{2}$ and where $+$ denotes the Minkowski sum. The {\it Newton polygon} of $N$ is the union of the compact edges of the boundary of $N$. Let $E$ be an edge of the Newton polygon of $N$. The {\it length} of $E$ (resp. {\it height} of $E$) is the length of the projection of $E$ onto the horizontal (respectively vertical axis). We denote by ${\rm inc}(E)$ the {\it inclination} of $E$, that is, the quotient of the length of $E$ by the height of $E$.

The {\it support} of any power series $f(x,y)=\sum_{ij}a_{ij}x^{i}y^{j}\in K[[x,y]]$ is 
$\hbox{\rm supp}(f):=\{(i,j)\in \mathbb N^{2}\;:\;a_{ij}\neq 0\}$. The Newton diagram of $f$ denoted by $\mathcal N(f)$  is the Newton diagram of its support. 
We can extend the notion of Newton diagram to any subset $S\subset \mathbb Q^{2}$ such that
there exists a positive integer $n$ with $n\cdot S:=\{n\cdot s\;:\;s\in S\}\subset \mathbb N^{2}$. 

A {\it supporting line} of $\mathcal N(f)$ is any line of the form $ax+by=c$ that intersects  $\mathcal N(f)$ and satisfies $ai+bj\geq c$ for any $(i,j)\in \mathcal N(f)$. Therefore a supporting line of $\mathcal N(f)$ either intersects this diagram at a vertex or along an edge of its Newton polygon. The inclination of the supporting line containing an edge $E$ of $\mathcal N(f)$ is, by definition, the inclination of $E$. Let $E$ be a compact edge of the Newton polygon of $f$. 
The \textit{Newton principal part} of $f$ with respect to $E$ is by definition $f_{\vert E}(x,y)=\sum_{(i,j)\in E\cap \hbox{\rm supp}(f)}a_{ij}x^{i}y^{j}$.

\begin{lemma}\label{ord}
Let $n$ be a positive integer and $g$ be a nonzero power series in $K[[x^{1/n},y]]$. Let $\nu$ be a positive rational number. Let $E$ be the intersection of the Newton diagram $\mathcal N(g)$ with its supporting line of inclination $\nu$  and let $(a,0)$ be the intersection point of this supporting line with the horizontal axis. Then for every $c\in K$ we have  $\ord g(x,cx^\nu) \geq a$, with equality if and only if $g_{\vert_E}(x,cx^\nu)\neq 0$.
\end{lemma}

\begin{proof}
Consider a nonzero power series $g(x,y)=\sum_{(i,j)\in (\frac{1}{n}\N)\times\N} a_{i,j}x^iy^j$. The set $D:=\{\, i+\nu j: (i,j)\in \hbox{\rm supp}(g) \,\}$ is well-ordered and its minimal element is $a$. For every $d\in D$ put $g_d(x,y)=\sum_{i+\nu j=d} a_{i,j}x^iy^j$. Then $g_d(x,cx^{\nu})=g_d(1,c)x^d$ and $g_{\vert_E}=g_a$. By equality $g=\sum_{d\in D}g_d$, we get $g(x,cx^{\nu})=g_a(x,cx^{\nu}) + \mbox{terms of orders bigger than $a$}$, which proves the lemma. 
\end{proof}

\section{Eggers decomposition of a power series }\label{sec:Ed}

In this section, $f\in K[[x,y]]$ is a power series such that $0 <\ord f(0,y)<\infty$. Let $f=\prod_{i=1}^r f_i$ be the factorization of $f$ into its irreducible factors and denote by $R$ the root of $\Theta(f)$.

Let $\Gamma=\prod_{j=1}^s \Gamma_j$ be the factorization into irreducible factors of some $\Gamma\in K[[x,y]]$ with $\Gamma(0,0)=0$. Given $P\in\Theta(f)$, denote
\begin{equation*}
\Gamma^{(P)}:=\prod \Gamma_j, 
\end{equation*}
where the product runs over $j\in\{1,\ldots, s\}$ such that the attaching point of $\Gamma_j$ to $\Theta(f)$ is $P$ (by convention the empty product is $1$).

\begin{definition}
 We say that the pair $(f,\Gamma$) satisfies condition $(E1)$ if 
\begin{itemize}
\item[(i)] $\Gamma=\prod_P \Gamma^{(P)}$, where the product runs over marked points $P$ of $\Theta(f)$ different from $R$, that is, the attaching points of irreducible factors of $\Gamma$ are marked points of $\Theta(f)$ different to $R$.
\item[(ii)] If $P$ is a leaf of $\Theta(f)$, then $i_0(\Gamma^{(P)},x)=i_0(f_P,x)-\de(P)$, otherwise
             $i_0(\Gamma^{(P)},x)=\sum_{P^*} \de(P^*) -\de(P)$, where the sum runs over direct successors $P^*$ of $P$.
\end{itemize}
\end{definition}

Note that if the pair $(f,\Gamma$) satisfies condition $(E1)$, $P$ is a leaf of $\Theta(f)$, and $f_P$ is irreducible then $i_0(\Gamma^{(P)},x)=0$.

\begin{definition}
We say that the pair $(f,\Gamma)$ satisfies condition $(E2)$ if for every $Q\in\Theta(f)$, $i_0(\Gamma_Q,x)=i_0(f_Q,x)-\de(Q)$.
\end{definition}

\begin{theorem}\label{thm:E1iffE2} 
The pair $(f,\Gamma)$ satisfies condition (E1) if and only if it satisfies condition (E2).
\end{theorem}

\begin{proof}
First we will prove that (E1) implies (E2). Assume that (E1) holds. According to item (i) of~(E1), 
\begin{equation}\label{eq:E1i}
    i_0(\Gamma,x)=\sum_P i_0(\Gamma^{(P)},x),
\end{equation}
where the sum runs over all marked points $P\in\Theta(f)$ different from $R$. 
Substituting the formulas from item (ii) of (E1) into~\eqref{eq:E1i} and canceling the terms $\de(P)$ of opposite signs we get 
\[ i_0(\Gamma,x) = -1 + \sum_P i_0(f_P,x)=-1+i_0(f,x), \]
where the sum runs over the leaves $P$ of $\Theta(f)$. 
Therefore, (E2) holds for $Q=R$ (note that $\Gamma=\Gamma_R$ and $f=f_R$).

For any $Q\in\Theta(f)$ different from $R$, by item (i) of (E1) we have $\Gamma_Q=\prod_P \Gamma^{(P)}$, 
where the product runs over marked points $P$ of $\Theta(f)$ such that $Q\preceq_x P$. Hence $i_0(\Gamma_Q,x)=\sum_{Q\preceq_x P} i_0(\Gamma^{(P)},x)$. Using the same argument as in the previous part of the proof we obtain  
\[ i_0(\Gamma_Q,x) = -\de(Q) + \sum_P i_0(f_P,x),\] 
where the sum runs over the leaves $P$ of $\Theta(f)$ such that $Q\preceq_x P$. Hence, $i_0(\Gamma_Q,x) =-\de(Q)+i_0(f_Q,x)$.

Next, we will prove that (E2) implies (E1). In what follows, we assume that (E2) holds.  Let us prove that item (i) of (E1) holds. Consider any two consecutive marked points $P'$ and $P$ of $\Theta(f)$. Then by condition (E2) the function 
$Q\to i_0(\Gamma_Q,x)$ is constant and equal to $i_0(f_Q,x)-\de(Q)=i_0(f_P,x)-\de(P)$ in the interval $(P',P]$. 
Therefore,  by Proposition~\ref{prop:step} there are no attaching points of the irreducible factors of $\Gamma$ in the interval $(P',P)$. 
If $P'=R$ then $f_R=f_P$ because $x$ is not a factor of $f$. Hence by (E2), we get $i_0(\Gamma_R,x)=i_0(f_R,x)-\de(R)=i_0(f_P,x)-\de(P)=i_0(\Gamma_P,x)$, 
so the function $Q\to i_0(\Gamma_Q,x)$ is constant in the interval $[R,P]$. It follows from Proposition~\ref{prop:step} that $R$ is not 
the attaching point to $\Theta(f)$ of any irreducible factors of $\Gamma$.

If $P$ is a leaf of $\Theta(f)$, then $\Gamma^{(P)}=\Gamma_P$ and the formula in item (ii) of (E1) holds.

Finally, let $P$ be a marked point of $\Theta(f)$ which is not a leaf. Let $P_{j_1}, \ldots, P_{j_t}$ be the marked points of $\Theta(f)$ that are direct successors of $P$. Since item (i) of (E1) holds, we deduce that $\Gamma_P=\Gamma^{(P)}\prod_{k=1}^t \Gamma_{P_{j_k}}$. Hence
    \[ i_0(\Gamma^{(P)},x)=i_0(\Gamma_P,x)-\sum_{k=1}^t i_0(\Gamma_{P_{j_k}},x)=i_0(f_P,x)-\de(P)-\sum_{k=1}^t\left(i_0(f_{P_{j_k}},x)-\de(P_{j_k})\right),\] 
where the last equality follows from the formula in (E2). Since $f_P=\prod_{k=1}^tf_{P_{j_k}}$, we obtain that $i_0(f_P,x)=\sum_{k=1}^ti_0(f_{P_{j_k}},x)$ and item (ii) of (E1) follows. 
\end{proof}

\begin{remark}\label{rem:computingd}
Let $P$ be any point of $\Theta(f)$ different from $R$. Consider any irreducible power series $h\in K[[x,y]]$ whose attaching point to $\Theta(f)$ is $P$. Then the logarithmic distance between $h$ and any irreducible factor $f_i$ of $f$ is determined by $\Theta(f)$ and $P$. Indeed, we have $d(f_i,h)=\ic(\langle x,f_i,h \rangle)$ and 
\begin{equation} \label{eq:cont}
d(f_i,h) =\left\{\begin{array}{ll} \ic(P) &\;\hbox{\rm if $f_i$ is above $P$}\\
\ic(\langle x,f_i,f_{i_0} \rangle)&\;\hbox{\rm otherwise,}
\end{array}
\right.
\end{equation}
where $f_{i_0}$ is any irreducible factor of $f$ above $P$. Notice that for any irreducible factor $\Gamma_j$ of $\Gamma^{(P)}$, the logarithmic distance $d(f_i,\Gamma_j)$ can be computed using \eqref{eq:cont}.
\end{remark}

\begin{definition}\label{def:Ed}
Let $f\in K[[x,y]]$ satisfy $0<\ord f(0,y)<\infty$ and let $\Gamma\in K[[x,y]]$ satisfy $\Gamma(0,0)=0$. We say that $\Gamma$ \textit{admits Eggers decomposition with respect to $f$} if one of the equivalent conditions $(E1)$ or $(E2)$ is satisfied by the pair $(f,\Gamma)$. In particular, if $1<\ord f(0,y)<\infty$ then we say that $\frac{\partial f}{\partial y}$ \textit{admits Eggers decomposition} if $\frac{\partial f}{\partial y}$ admits Eggers decomposition with respect to $f$.
\end{definition}

The following corollary compiles the information on the factors of the Eggers decomposition of a power series. Its proof follows from  Theorem \ref{thm:E1iffE2} and Remark \ref{rem:computingd}. Note that this information only depends on $\Theta(f)$. 

\begin{corollary}\label{coro:Edecg}
Let $f\in K[[x,y]]$ be such that $0<\ord f(0,y)<\infty$ and let $\Gamma\in K[[x,y]]$ with $\Gamma(0,0)=0$. If $\Gamma$ admits Eggers decomposition with respect to $f$, then $\Gamma=\prod_{P}\Gamma^{{(P)}}$ where the product runs over marked points $P$ of $\Theta(f)$ different from $R$, and we have 
\begin{enumerate}
\item $i_0(\Gamma^{(P)},x)=\left\{\begin{array}{ll}
i_0(f_P,x)-\de(P) & \hbox{\rm if $P$ is a leaf}\\
 & \\
\sum_{P^*\in {\mathcal S}_{P}} \de(P^*) -\de(P) & \hbox{\rm otherwise},
\end{array}
\right.$

\noindent where ${\mathcal S}_{P}$ denotes the set of direct successors of $P$ in $\Theta (f)$.

\item  If $f_{i}$ is an irreducible factor of $f$, then 
\[
\frac{i_0(\Gamma^{(P)},f_{i})}{i_0(\Gamma^{(P)},x)}=\left\{\begin{array}{ll}
 i_{0}(f_{i},x)\ic(P)&  \hbox{\rm if $f_{i}$ is above $P$}\\
  & \\
  i_{0}(f_{i},x)\ic(\langle x,f_i,f_{i_0} \rangle)&  \hbox{\rm otherwise}\\
\end{array}
\right. 
\]
\noindent where $f_{i_0}$ is any irreducible factor of $f$ above $P$.
\end{enumerate}
\end{corollary}

The following two results will be useful in Section \ref{sec:5}. Let us keep the notation introduced in this section. Let $Q$ be a point of $\Theta(f)$ and let $h\in K[[x,y]]$ be irreducible. We set
\[ d(h,Q) := \ic(\langle x, P_h , Q \rangle),\]
where $P_h$ denotes the attaching point of $h$ to $\Theta(f)$.

\begin{theorem}\label{thm:1}
If the pair $(f,\Gamma)$ satisfies condition (E2), then for any point $Q$ of $\Theta(f)$ different from the leaves we have
\begin{equation}\label{E:ex}  
\sum_{i=1}^r d(f_i,Q)i_0(f_i,x)-\sum_{j=1}^s d(\Gamma_j,Q)i_0(\Gamma_j,x)=\ex(Q).
\end{equation}
\end{theorem}

\begin{proof}
We will show that formula~\eqref{E:ex} prolongs along the Eggers-Wall tree of $f$. If $Q$ is the root of $\Theta(f)$, then both sides of \eqref{E:ex} are equal to zero. 
Next consider two consecutive marked points $P$ and $P'$ of $\Theta(f)$ and suppose that \eqref{E:ex} holds for $P$. Let $Q$ be any point of $[P,P']$ which is not a leaf. Put $S:=\sum_{i=1}^r d(f_i,Q)i_0(f_i,x)-\sum_{j=1}^s d(\Gamma_j,Q)i_0(\Gamma_j,x)-\ex(P)$. Since  \eqref{E:ex} holds for $P$ we get
\begin{equation}\label{eq:xxx}
S=\sum_{i=1}^r (d(f_i,Q)-d(f_i,P))i_0(f_i,x)-
\sum_{j=1}^s (d(\Gamma_j,Q)-d(\Gamma_j,P))i_0(\Gamma_j,x).
\end{equation}

\noindent If $Q\not\in[R,P_{f_i}]$ then $d(f_i,Q)=d(f_i,P)$. By Theorem \ref{thm:E1iffE2}, the pair $(f,\Gamma)$ satisfies item $(i)$ of $(E1)$, hence the attaching point $P_{\Gamma_j}$ of $\Gamma_j$ to $\Theta(f)$ is a marked point different from $R$, for $1\leq j\leq s$. Moreover if $Q\not\in[R,P_{\Gamma_j}]$ then $d(\Gamma_j,Q)=d(\Gamma_j,P)$, so \eqref{eq:xxx} becomes

\begin{equation}\label{eq:yyy}
S=\sum_{Q\preceq_x P_{f_i}} (d(f_i,Q)-d(f_i,P))i_0(f_i,x)-\sum_{Q\preceq_x P_{\Gamma_j}} (d(\Gamma_j,Q)-d(\Gamma_j,P))i_0(\Gamma_j,x).
\end{equation}

\noindent Observe that by definition, if $Q\in[R,P_{f_i}]$ then $d(f_i,P)=\ic(P)<\ic(Q)=d(f_i,Q)$, and the same is true if we replace $f_i$ with $\Gamma_j$, therefore from \eqref{eq:yyy} we get

\begin{eqnarray*}
S&=&\sum_{Q\preceq_x P_{f_i}} (\ic(Q)-\ic(P))i_0(f_i,x)-
\sum_{Q\preceq_x P_{\Gamma_j}} (\ic(Q)-\ic(P))i_0(\Gamma_j,x)\\
&=&(\ic(Q)-\ic(P))\left(\sum_{Q\preceq_x P_{f_i}} i_0(f_i,x) - \sum_{Q\preceq_x P_{\Gamma_j}} i_0(\Gamma_j,x)\right)\\
&=&(\ic(Q)-\ic(P))(i_0(f_Q,x)-i_0(\Gamma_Q,x))=(\ic(Q)-\ic(P))\de(Q)=\ex(Q)-\ex(P),
\end{eqnarray*}

\noindent where the fourth equality holds because $(f,\Gamma)$ satisfies condition $(E2)$ and the fifth equality holds by the definition of the exponent function (see \eqref{def:ex}).
We conclude that equality \eqref{E:ex} holds for $Q$.
\end{proof}

\begin{corollary}\label{co:pcondition}
    Assume that the pair $(f,\Gamma)$ satisfies condition (E2). Let $g\in K[[x]][y]$ be the minimal polynomial of a Puiseux series $\alpha(x)$ of positive order such that $f(x,\alpha(x))\neq 0$, and let $Q$ be the attaching point of $g$ to $\Theta(\Gamma f)$. If $Q$ is a point of $\Theta(f)$, then $\ord f(x,\alpha(x))-\ord \Gamma(x,\alpha(x)) = \ex(Q).$
\end{corollary}

\begin{proof}
Since $\ord f_i(x,\alpha(x))=d(f_i,g)i_0(f_i,x)=d(f_i,Q)i_0(f_i,x)$ for $1\leq i\leq r$ 
and $\ord \Gamma_j(x,\alpha(x))=d(\Gamma_j,g)i_0(\Gamma_j,x)=d(\Gamma_j,Q)i_0(\Gamma_j,x)$ for $1\leq j\leq s$, we have 
\begin{eqnarray*}
\ord f(x,\alpha(x))-\ord \Gamma(x,\alpha(x)) &=&
\sum_{i=1}^r \ord f_i(x,\alpha(x)) -\sum_{j=1}^s \ord \Gamma_j(x,\alpha(x))  \\
&=& \sum_{i=1}^r d(f_i,Q)i_0(f_i,x)-\sum_{j=1}^s d(\Gamma_j,Q)i_0(\Gamma_j,x)=\ex(Q),
\end{eqnarray*}
where the last equality follows from Theorem \ref{thm:1}.
\end{proof}

The aim of the following two sections will be to show that $\frac{\partial f}{\partial y}$ admits Eggers decomposition if and only if $f$ verifies the Eggers condition. 

\section{Eggers condition is sufficient for Eggers decomposition}
\label{sec:sufficient}

We will prove in this section  that if $f\in K[[x,y]]$, with $1<\ord f(0,y)<\infty$, verifies the Eggers condition then $\frac{\partial f}{\partial y}$ admits Eggers decomposition (see Corollary \ref{Eggersfact}). This will be a consequence of the first theorem we present in this paper (see Theorem \ref{mainJ}), which under certain assumptions related to the Eggers-Wall tree of $f$, guarantees a decomposition of $\Gamma:=\frac{\partial f}{\partial y}$ as products of series $\Gamma^{(P)}$ that are not necessarily irreducible, where the intersection multiplicities of $\Gamma^{(P)}$ with $x$ and any irreducible factor of $f$ (and so with $f$) only depend on the functions defined on the Eggers-Wall tree of $f$. We also present in this section several examples illustrating this theorem. Let us start with some lemmas necessary for our objective.

\begin{lemma}\label{EdgeE}
Let $f\in K[[x,y]]$ be an irreducible power series such that $0<\ord f(0,y)<\infty$ and let $\alpha(x)$ be a Puiseux series of positive order such that $f(x,\alpha(x))\neq 0$. 
Let $g\in K[[x]][y]$ be the minimal polynomial of $\alpha(x)$. Assume that $\deg(g)\notequiv 0 \pmod p$. 
Let $P$ be the ramification point of the Eggers tree $\Theta(fg)$. 
Set $\tilde f(x,y)=f(x,y+\alpha(x))$. Then, the Newton polygon of $\tilde f$ has at least one compact edge. Moreover, if $E$ is the compact edge of $\mathcal N(\tilde f)$ with maximal inclination, then:
\begin{enumerate}
\item The inclination of $E$ equals $\ex(P)$.
\item The endpoints of $E$ are $(a,i_0(f,x)/\de(P))$ 
      and $(d(f,g) i_0(f,x),0)$, with $a\geq0$.
\end{enumerate}
\end{lemma}

\begin{proof}
The Newton polygon of $\tilde f$ has at least one compact edge because $0<\ord\tilde f(x,0)<\infty$ and $0<\ord \tilde f(0,y)<\infty$.
The point $B:=(d(f,g)i_0(f,x),0)$ is the endpoint of $E$ which lies on the horizontal axis. 
This is true By equality $\ord\tilde f(x,0)=\ord f(x,\alpha(x))=i_0(f,g)/i_0(g,x)=d(f,g)i_0(f,x)$, Denote by $A$ the second endpoint of $E$ and by ${\rm inc}(E)$ the inclination of $E$.

Let $J$ be a set of quotients $m/n$ such that $m,n$ are positive integers and $\gcd(m,n)=\gcd(n,p\deg g)=1$  and  there is no compact edge of $\mathcal N(\tilde f)$ of inclination $m/n$. The set $J$ is dense in the set of positive rational numbers. Since the elements of the support of $\alpha(x)$ are of the form $i/\deg g$ and $\gcd(m,n)=1$, we get that $J$ is disjoint with the support of $\alpha(x)$.

For every $\iota\in J$, we consider the Puiseux series $\gamma_{\iota}(x)=\alpha(x)+x^{\iota}$ and the minimal polynomial $g_{\iota}\in K[[x]][y]$ of $\gamma_{\iota}(x)$. It follows from the choice of $\iota$ that $\deg(g_\iota)\notequiv 0 \pmod p$. We have 
\begin{equation}\label{dfgiota}
    \ord \tilde f(x,x^{\iota})=\ord f(x,\gamma_{\iota}(x))=d(f,g_{\iota}) i_0(f,x).
\end{equation} 
Since $\iota$ does not belong to the support of $\alpha(x)$, we have 
\[ k(g,g_\iota)=\max\set{\ord(\gamma(x)-\alpha(x))\;:\;\gamma(x)\text{ Puiseux root of } g_\iota}=\iota.
\]
Hence, if we denote by $Q_\iota$ the point $\langle x, g, g_\iota \rangle$ of the tree $\Theta(fgg_\iota)$, then by Property \ref{prop:1}, $\ex(Q_\iota)=\iota$.

In order to show that ${\rm inc}(E)=\ex(P)$ it is enough to exclude the cases ${\rm inc}(E)<\ex(P)$ and ${\rm inc}(E)>\ex(P)$. In what follows, the point $P$ is viewed as the point $\langle x, f, g \rangle$ of the Eggers-Wall tree $\Theta(fgg_\iota)$. Figure \ref{fig:1} contains a drawing of $\Theta(fgg_\iota)$ in each case.

\begin{figure}
    \begin{center}
\begin{tikzpicture}[scale=1]
\begin{scope}[shift={(0,0)}]

    \draw [-, color=black,  line width=1pt](0,0) -- (0,2);
    \draw [->, color=black,  line width=1pt](0,2) -- (1,3.3);
    \draw [->, color=black,  line width=1pt](0,2) -- (-1,3.3);
    \draw [->, color=black,  line width=1pt](0,1) -- (0,0);
    \draw [->, color=black,  line width=1pt](0,1) -- (1,2.2); 
     
    \node [below, color=black] at (0,0) {$x$};
    \node [above, color=black] at (-1,3.3) {$f$};
    \node [above, color=black] at (1,3.3) {$g$};
    \node [right, color=black] at (1,2.2) {$g_{\iota}$};
              
    \node[draw,circle, inner sep=1.5pt,color=black, fill=black] at (0,1){};
    \node[draw,circle, inner sep=1.5pt,color=black, fill=black] at (0,2){};
    \node [left, color=black] at (0,1) {$Q_\iota$}; 
    \node [left, color=black] at (0,2) {$P$}; 
    
    \end{scope}
      
    \begin{scope}[shift={(5,0)}]

    \draw [-, color=black,  line width=1pt](0,0) -- (0,2);
    \draw [->, color=black,  line width=1pt](0,2) -- (1,3.3);
    \draw [->, color=black,  line width=1pt](0,2) -- (-1,3.3);
    \draw [->, color=black,  line width=1pt](0,1) -- (0,0);
    \draw [->, color=black,  line width=1pt](0.5,5.3/2) -- (1.5,3);
     
    \node [below, color=black] at (0,0) {$x$};
    \node [above, color=black] at (-1,3.3) {$f$};
    \node [above, color=black] at (1,3.3) {$g$};
    \node [right, color=black] at (1.5,3) {$g_{\iota}$};

    \node[draw,circle, inner sep=1.5pt,color=black, fill=black] at (0.5,5.3/2){};
    \node[draw,circle, inner sep=1.5pt,color=black, fill=black] at (0,2){};
      
    \node [left, color=black] at (0.4,5.3/2) {$Q_\iota$}; 
    \node [left, color=black] at (0,1.9) {$P$}; 
    \end{scope}
  \end{tikzpicture}
\end{center}
 \caption{Both situations on the proof of Lemma \ref{EdgeE}}
\label{fig:1}
   \end{figure}
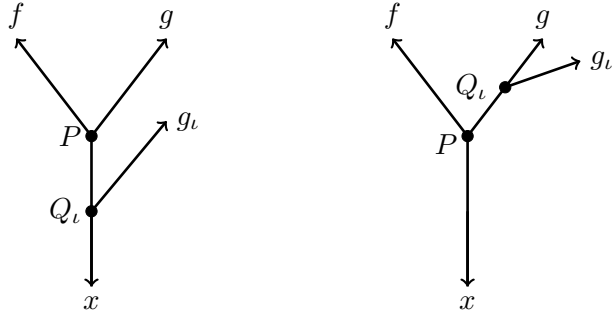

First suppose that ${\rm inc}(E)<\ex(P)$. Then there exists $\iota\in J$ such that ${\rm inc}(E)<\iota<\ex(P)$. Since $\ex(Q_\iota)=\iota<\ex(P)$, we have also $\ic(Q_\iota)<\ic(P)$. By inequality ${\rm inc}(E)<\iota$, the supporting line of $\mathcal N(\tilde f)$ of inclination $\iota$ intersects $\mathcal N(\tilde f)$ only at $B$. Hence, from Lemma~\ref{ord}  and equalities \eqref{dfgiota}, it follows that $d(f,g) i_0(f,x)= \ord\tilde f(x,x^\iota) = d(f,g_\iota) i_0(f,x)$, which gives $\ic(P)=d(f,g)=d(f,g_\iota)=\ic(Q_\iota)$. This is contradiction.

Now suppose that ${\rm inc}(E)>\ex(P)$. Take $\iota\in J$ such that ${\rm inc}(E)>\iota>\ex(P)$. On the one hand, we have $\ex(Q_\iota)=\iota>\ex(P)$ and consequently $\ic(P)=d(f,g)=d(f,g_\iota)$. On the other hand, since ${\rm inc}(E)>\iota$, the supporting line of $\mathcal N(\tilde f)$ of inclination $\iota$ intersects $\mathcal N(\tilde f)$ at a vertex different from~$B$. By Lemma~\ref{ord} and equalities \eqref{dfgiota}, we have $d(f,g_{\iota})i_0(f,x) = \ord\tilde f(x,x^\iota) < d(f,g)i_0(f,x)$ This is a contradiction. 

Finally, let us calculate  the second coordinate $h$ of the point $A$. Denote by $P'$ the marked point in $\Theta(fg)$ such that $P'\prec_x P$ and there are no marked points in $(P',P)$. 
Take $\iota\in J$ ($\ex(P')<\iota<\ex(P)$) such that  the supporting line $\ell$ of $\mathcal N(\tilde f)$ of inclination $\iota$ intersects $\mathcal N(\tilde f)$ only at $A$. 
Let $B'$ be the point where $\ell$ intersects the horizontal axis. Consider the triangles $ABB'$ and $Abb'$ (see Figure \ref{fig:2}). 

\begin{figure}[ht]
    \begin{center}
\begin{tikzpicture}[scale=1]
 \draw [->, color=black,  line width=1pt](0,0) -- (0,2);
 \draw [->, color=black,  line width=1pt](0,0) -- (4,0);
 \draw [-, color=black,  line width=1pt](1,1) -- (3,0);
 \draw [-, dashed, color=black,  line width=1pt](0.5,2) -- (2,-1);
 \draw [-, dashed, color=black,  line width=1pt](1,2) -- (1,0);
 \draw [-, dashed, color=black,  line width=1pt](1,1) -- (-1,1);
 \draw [-, dashed, color=black,  line width=1pt](1,0.5) -- (2,0.5);

\node[draw,circle, inner sep=1.5pt,color=black, fill=black] at (1,1){};
\node[draw,circle, inner sep=1.5pt,color=black, fill=black] at (3,0){};
\node[draw,circle, inner sep=1.5pt,color=black, fill=black] at (1.5,0){};
\node[draw,circle, inner sep=1.5pt,color=black, fill=black] at (2,0.5){};
\node[draw,circle, inner sep=1.5pt,color=black, fill=black] at (5/4,0.5){};
 
\node [above, color=black] at (1.1,1) {$A$};
\node [below, color=black] at (3,0) {$B$};
\node [below,  color=black] at (1.4,-0.03) {$B'$};
\node [left, color=black] at (-0.2,0.5) {$h$};
\node [below, color=black] at (2,0.5) {$b$};
\node [above,  color=black] at (2.5,1.55) {$b'$};
\node [right, below, color=black] at (0.7,0.95) {$1$};

\draw[<->] (0.9,0.9) to [bend right] (0.9,0.5);
\draw[<->] (-0.1,1) to [bend right] (-0.1,0);

\draw[->][thick, color=black](2.5,1.6) .. controls (2,1) ..(5/4,0.6);
\end{tikzpicture}
\end{center}
\caption{The edge $E$ and the triangles $ABB'$ and $Abb'$ of the proof of Lemma \ref{EdgeE}.}
\label{fig:2}
\end{figure}
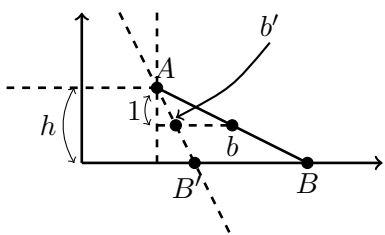

Since $\ex(P')<\ex(Q_\iota)=\iota<{\rm inc}(E)=\ex(P)$, 
we get $P'\prec_x Q_\iota=\langle x, g, g_\iota \rangle=\langle x, f, g_\iota \rangle\prec_x P$ in the tree $\Theta(fgg_\iota)$ and 
\[\vert\overline{BB'}\vert=\ord\tilde f(x,0)-\ord\tilde f(x,x^\iota)=(\ic(P)-\ic(Q_\iota))i_0(f,x)=\frac{\ex(P)-\ex(Q_\iota)}{\de(P)}i_0(f,x).
\]
To end the proof, we use an argument that relies on similarity of triangles. We obtain that  $h=\frac{\vert\overline{B'B}\vert}{\vert\overline{b'b}\vert}=\frac{\ex(P)-\ex(Q_\iota)}{\ex(P)-\ex(Q_\iota)}\frac{i_0(f,x)}{\de(P)}$, and the lemma follows.
\end{proof}

Although the following result will not be useful until Section \ref{sec:5}, we prove it here because it is a consequence of Lemma \ref{EdgeE}. We keep the notation of that lemma.

\begin{corollary}\label{C2} 
Under the assumptions of Lemma \ref{EdgeE} there exists $c\neq 0$ such that for the minimal polynomial $g_{\beta}\in K[[x]][y]$ of $\beta(x)=\alpha(x)+cx^{\ex(P)}$ we have $d(f,g_{\beta})>\ic(P)$.
\end{corollary}

\begin{proof}
By Lemma \ref{EdgeE} applied to $f$, $\alpha(x)$, and $g$, the Newton polygon of $\tilde f$ has a compact edge $E$ of inclination $\ex(P)$ and whose intersection with the horizontal axis is the point $(d(f,g)i_0(f,x),0)$. Since $\tilde f_{\vert E}(1,y)\in K[y]$ is not a monomial, we can choose $c\in K^*$ such that $\tilde f_{\vert E}(1,c)=0$. The Newton principal part of $\tilde f$ with respect to $E$  satisfies $\tilde f_{\vert E}(x,cx^{\ex(P)})=\tilde f_{\vert E}(1,c)x^{d(f,g)i_0(f,x)}$, therefore $\tilde f_{\vert E}(x,cx^{\ex(P)})=0$ and it follows from Lemma \ref{ord} that $\ord \tilde f(x,cx^{\ex(P)})>d(f,g)i_0(f,x)$. We have  $d(f,g_{\beta})i_0(f,x)=i_0(f,g_{\beta})/i_0(g_{\beta},x)=\ord f(x,\beta(x))=\ord \tilde f(x,cx^{\ex(P)})$ and $\ic(P)=d(f,g)$, so this ends the proof.
\end{proof}

\begin{lemma}\label{edgeE23}
Let $f=\prod_{i=1}^r f_i$ be the factorization of $f\in K[[x,y]]$ into irreducible factors and assume that $0<\ord f(0,y)<\infty$. Let $[P',P]$ be a segment of $\Theta(f)$ with no marked points in $(P',P)$ and such that $\de(P)\notequiv 0 \pmod p$. Let $\iota=m/n$ be a rational number such that ${\rm gcd}(m,n)={\rm gcd}(p\cdot\de(P),n)=1$ and $e(P')<\iota<e(P)$. Then:
\begin{enumerate}
    \item\label{edgeE23:s1} There exists a Puiseux series $\alpha(x)=\sum_{i/\de(P)<\iota} a_ix^{i/\de(P)}-x^\iota$ such that for the minimal polynomial $g\in K[[x]][y]$ of $\alpha(x)$, the attaching point $Q$ of $g$ to $\Theta(f)$ belongs to $(P',P)$ and $e(Q)=\iota$.
    \item\label{edgeE23:s2} Set $\tilde f(x,y)=f(x,y+\alpha(x))$. Then, the Newton polygon of $\tilde f$ has at least one compact edge. Moreover, if $E$ is the compact edge of maximal inclination of $\mathcal N(\tilde f)$, then:
    \begin{enumerate}[label=\normalfont(\alph*)]
    \item\label{edgeE23:s21} The inclination of $E$ is $\iota$.
    \item\label{edgeE23:s22} The endpoints of $E$ are $(a,i_0(f_P,x)/\de(P))$ and $(\sum_{i=1}^r d(f_i,g)i_0(f_i,x),0)$, with $a\geq0$.
    \item\label{edgeE23:s23} The Newton principal part of $\tilde f$ with respect to $E$ is of the form $\tilde f_{\vert E}=\varepsilon x^a(y-x^\iota)^{i_0(f_P,x)/\de(P)}$, with $\varepsilon\in K^*$.
    \end{enumerate}
    \end{enumerate}
\end{lemma}

\begin{proof}
Since $P$ is different from the root of the tree $\Theta(f)$, there exists a Weierstrass polynomial $h\in K[[x]][y]$ such that $h$ is above $P$ and $\deg h=\de(P)$ 
(it is enough to choose a convenient key polynomial of a factor of $f$ which is above $P$). Let $\gamma(x)=\sum_{i>0} a_ix^{i/\de(P)}$ be a Puiseux root of $h$ and set $\alpha(x)=\sum_{i/\de(P)<\iota} a_ix^{i/\de(P)}-x^\iota$. Let $g\in K[[x]][y]$ be the minimal polynomial of $\alpha(x)$. Observe that $\deg(g)=n\cdot\de(P)\notequiv 0 \pmod p$ and $\iota\notin{\rm supp}(\gamma(x))$. Hence, by Property \ref{prop:1}, we have 
\[\iota=\max\set{\ord(\beta(x)-\gamma(x))\;:\;\beta(x)\text{ Puiseux root of }g}=\ex(Q),\] where $Q$ is the ramification point of the Eggers tree $\Theta(hg)$. The point $Q$ can be viewed in the tree $\Theta(f g)$ as the point $\langle x, f_{i}, g \rangle$ for any factor $f_i$ of $f$ which is above $P$, and $Q\in(P',P)$ because $\ex(P')<\iota=\ex(Q)<\ex(P)$.  

Next we prove statement \eqref{edgeE23:s2}. The fact that $\mathcal N(\tilde f)$ has at least one compact edge follows from the same arguments given in the proof of Lemma \ref{EdgeE}.

Let us prove statements \ref{edgeE23:s21} and \ref{edgeE23:s22}. For every $i\in\set{1,\ldots,r}$, set $\tilde f_i(x,y)=f_i(x,y+\alpha(x))$. Since $\deg(g)\notequiv 0 \pmod p$ and $f_i(x,\alpha(x))\neq0$, we are under the assumptions of Lemma \ref{EdgeE}. Let $E_i$ be the compact edge of maximal inclination of $\mathcal N(\tilde f_i)$. It follows that ${\rm inc}(E_i)=\ex(\langle x,f_i,g\rangle)$, and the endpoints of $E_i$ are $A_i=(a_i,i_0(f_i,x)/\de(\langle x,f_i,g\rangle))$ and $B_i=(d(f_i,g)i_0(f_i,x),0)$. If $f_i$ is above $P$, then $\langle x,f_i,g\rangle=Q$ and thus ${\rm inc}(E_i)=\iota$ and $\de(\langle x,f_i,g\rangle)=\de(Q)=\de(P)$. Otherwise, $\langle x,f_i,g\rangle \prec_x Q$ and we get ${\rm inc}(E_i)<\iota$. Since the Newton diagram of $\tilde f$ is the Minkowski sum $\sum_{i=1}^r\mathcal N(\tilde f_i)$, we deduce that ${\rm inc}(E)=\iota$ and $(\sum_{i=1}^r d(f_i,g)i_0(f_i,x),0)$ is the endpoint of $E$ living on the horizontal axis. Moreover, if we denote by $A$ the second endpoint of $E$, then the second coordinate of $A$ is the sum of the second coordinates of endpoints $A_i$ where the sum runs on $i$ such that  $E_i$ have inclination $\iota$. This ends the proof of statements \ref{edgeE23:s21} and \ref{edgeE23:s22}.

Let us prove statement \ref{edgeE23:s23}. Since the inclination of $E$ is $\iota$ and $E$ intersects the horizontal axis, the Newton principal part $\tilde f_{\vert E}$ has the form  $\tilde f_{\vert E}=\varepsilon x^a \prod_{i=1}^{i_0(f_P,x)/\de(P)}(y-c_ix^\iota)$, with $\varepsilon, c_i\in K^*$. It is enough to prove that $c_i=1$ for any $i$. Take any $c\in K\backslash \{1\}$ and consider the minimal polynomial $g_c\in K[[x]][y]$ of $\alpha_c(x):=\alpha(x)+cx^\iota$. Since the attaching points of $g$ and $g_c$ to $\Theta(f)$ are both equal to $Q$, we get that $d(f_i,g_c)=d(f_i,g)$ for all $i\in\set{1,\ldots,r}$. Therefore, 
\[
\ord \tilde f(x,0)=\sum_{i=1}^rd(f_i,g)i_0(f_i,x)=\sum_{i=1}^rd(f_i,g_c)i_0(f_i,x)=\ord f(x,\alpha_c(x))=\ord \tilde f(x,cx^\iota).
\]
By Lemma \ref{ord} we have $\tilde f_{\vert E}(x,cx^\iota)\neq 0$. Hence, $c\neq c_i$ for every $i\in \{1, \ldots, i_0(f_P,x)/\de(P)\}$ and the lemma follows.
\end{proof}

\begin{theorem}\label{mainJ}
Let $f\in K[[x,y]]$ be a power series such that $1<\ord f(0,y)<\infty$ and denote $\Gamma=\frac{\partial f}{\partial y}$. Let $P'$ and $P$ be two consecutive marked points of $\Theta(f)$ such that $f$ verifies the Eggers condition at P. Then, for every $Q$ in the interval $(P',P]$, we have $i_0(\Gamma_Q,x)=i_0(f_Q,x)-\de(Q)$.
\end{theorem}

\begin{proof}
According to Proposition~\ref{prop:step}, the function $\phi$ that sends $Q$ to $i_0(\Gamma_Q,x)$ is a weakly decreasing step function which is continuous from the left at every point of $(P',P]$. Assume that there exists a dense subset $S$ of $(P',P]$ such that $\phi(Q)=i_0(f_Q,x)-\de(Q)$ for all $Q\in S$. Then the function $\phi$ is constant on $S$, because $f_Q=f_P$ and $\de(Q)=\de(P)$ for all $Q\in(P',P]$. This implies that $\phi$ is constant on $(P',P]$ and $i_0(\Gamma_Q,x)=i_0(f_P,x)-\de(P)=i_0(f_Q,x)-\de(Q)$ for every $Q$ in the interval $(P',P]$, so this ends the proof. In what follows, we prove that such a subset $S$ of $(P',P]$ exists.

Let $S$ be the subset of $(P',P]$ consisting of $Q\in(P',P)$ such that $\ex(Q)=\iota$, where $\iota=m/n$ is a rational number such that ${\rm gcd}(m,n)={\rm gcd}(p\cdot\de(P),n)=1$, and $Q$ is not the attaching point to $\Theta(f)$ of any irreducible factor of $\Gamma$. Observe that $S$ is a dense subset of $(P',P]$. Take $Q\in S$. We will show that $i_0(\Gamma_Q,x)=i_0(f_Q,x)-\de(Q)$.

Let $\Gamma=\frac{\partial f}{\partial y}=\prod_{j=1}^s \Gamma_j$ be the factorization into irreducible factors of $\Gamma$. By Properties \ref{prop:pE}, since $f$ verifies Eggers condition at $P$, it also verifies $\de$--condition at $P$, that is, $\de(P)\not\equiv0 \pmod p$.

Consider the Puiseux series $\alpha(x)$ and the series $\tilde f$ as in Lemma \ref{edgeE23} and denote by $g\in K[[x]][y]$ the minimal polynomial of $\alpha(x)$.   
Let $\tilde \Gamma(x,y):=\frac{\partial f}{\partial y}(x,y+\alpha(x))=\frac{\partial \tilde f}{\partial y}(x,y)$. 
If $E$ is the compact edge of $\mathcal N(\tilde f)$ of maximal inclination, then by Lemma \ref{edgeE23}--\ref{edgeE23:s23} we get $\frac{\partial }{\partial y}(\tilde f_{\vert E})=\varepsilon x^a k(y-x^\iota)^{k-1}$, with $\varepsilon\in K^*$ and $k=i_0(f_P,x)/\de(P)$. 
By the assumption $i_0(f_P,x)\notequiv 0 \pmod p$, we get $k\notequiv0 \pmod p$. Therefore, if $k>1$, 
then $\mathcal N(\tilde \Gamma)$ has a compact edge $\tilde E$ with ${\rm inc}(\tilde E)=\iota$ which intersects the horizontal axis, and the height of $\tilde E$ is equal to $k-1$. 
If $k=1$, then $\mathcal N(\tilde \Gamma)$ has no compact edge of inclination $\iota$.

We can write $\Gamma=x^d\Gamma'$ with $d\geq0$ where $\Gamma'$ is not divisible by $x$. We have $0<\ord\Gamma'(0,y)<\infty$ because $\ord f(0,y)>1$. Set $\tilde\Gamma'(x,y):=\Gamma'(x,y+\alpha(x))$. Observe that $\mathcal N(\tilde\Gamma)$ is the translation of $\mathcal N(\tilde\Gamma')$ by the horizontal vector $(d,0)$.

We assume first that there exists an irreducible factor $\Gamma_j$ of $\Gamma'$ which is above $Q$ in $\Theta(f)$. Since $Q$ is not the attaching point to $\Theta(f)$ of any irreducible factor of $\Gamma'$, we can choose $P_1',P_1\in(P',P)\subset\Theta(f)$ such that $P_1'\prec_x Q\prec_x P_1$ and the attaching points to $\Theta(f)$ of the irreducible factors of $\Gamma'$ do not belong to $[P_1',P_1]$. The segment $[P_1',P_1]$ can also be seen in the Eggers-Wall tree of $\Gamma'$ and the value of the index function (resp. of the exponent function) on this tree at any point $T\in[P_1',P_1]$ coincides with $\de(P)$ (resp. with $\ex(T)$) seen in $\Theta(f)$. Moreover, note that there are no marked points of $\Theta(\Gamma')$ in $(P_1',P_1)\subset\Theta(\Gamma')$. Hence Lemma~\ref{edgeE23} can be applied to the power series $\Gamma'$, the segment $[P_1',P_1]$, and the Puiseux series $\alpha(x)$. By statement \eqref{edgeE23:s2} of that lemma, we deduce that $k>1$ and $k-1=i_0(\Gamma_{P_1}',x)/\de(P_1)$. Since $\Gamma_{P_1}'=\Gamma_Q'=\Gamma_Q$, $f_Q=f_P$, $\de(P_1)=\de(Q)=\de(P)$, and $k=i_0(f_P,x)/\de(P)$, we obtain $i_0(\Gamma_Q,x)=i_0(f_Q,x)-\de(Q)$.

Suppose now that no irreducible factor of $\Gamma'$ is above $Q$ in $\Theta(f)$. For every irreducible factor $\Gamma_j$ of $\Gamma'$, set $\tilde \Gamma_j(x,y):=\Gamma_j(x,y+\alpha(x))$. Since $\deg(g)\notequiv 0 \pmod p$ and $\Gamma_j(x,\alpha(x))\neq0$, we are under the assumptions of Lemma \ref{EdgeE}. Therefore, ${\rm inc}(\tilde E_j)=\ex(\langle x,\Gamma_j,g\rangle)<\ex(Q)=\iota$, where $\tilde E_j$ is the compact edge of maximal inclination of $\mathcal N(\tilde \Gamma_j)$. Since the Newton diagram of a product is the Minkowski sum of the Newton diagrams of the factors, we obtain that 
$\mathcal N(\tilde \Gamma')$ has no compact edge of inclination $\iota$. It follows that $k=1$, that is, $i_0(f_P,x)-\de(P)=0$. Since $\Gamma_Q=\Gamma'_Q=1$, we obtain $i_0(f_Q,x)-\de(Q)=i_0(f_P,x)-\de(P)=0=i_0(\Gamma_Q,x)$, so the equation also holds in this case. 
\end{proof}

The following examples illustrate Theorem \ref{mainJ}.

\begin{example}\label{ex:p&p+1}
Let $K$ be an algebraically closed field of characteristic $k$. Consider $f(x,y)=y(y^{p}+x^{q})\in K[[x,y]]$. Here $p<q$ are coprime numbers. The Eggers-Wall tree of $f$ is drawn in Figure \ref{fig:ex}. In this figure $Q$ (respectively $T$)  is the leaf corresponding to $y$ (respectively to $y^{p}+x^{q}$). Notice that  $f_{P}=f$, $f_{Q}=y$, $f_{T}=y^{p}+x^{q}$, $\de(P)=\de(Q)=1$, and $\de(T)=p$. Therefore, $i_0(f_P,x)-\de(P)=p$ and $i_0(f_Q,x)-\de(Q)=i_0(f_T,x)-\de(T)=0$. We distinguish three cases:

\begin{enumerate}
\item If $k\not\in\{p,p+1\}$, then $\Gamma:=\frac{\partial f}{\partial y }=(p+1)y^{p}+x^{q}$ whose attaching point to $\Theta(f)$ is $P$. We get $\Gamma_{Q}=\Gamma_{T}=1$ and $\Gamma_{P}=\Gamma$. Observe that $f$ satisfies the Eggers condition at the points $P$, $Q$, and $T$.

\item If $k=p$, then $\frac{\partial f}{\partial y }$ is the branch of $f$ corresponding to the leaf $T$. Hence $\Gamma_{Q}=1$ and $\Gamma_{P}=\Gamma_{T}=y^{p}+x^{q}$. The power series $f$ does not satisfy the Eggers condition at the point $T$ and $i_0(\Gamma_T,x)\neq i_0(f_T,x)-\de(T)$. Nevertheless $f$ verifies the Eggers condition at both $P$ and $Q$.

\item Suppose that $k=p+1$. We get $\Gamma:=\frac{\partial f}{\partial y }=x^{q}$. Hence $\Gamma_X=1$ for any $X\in \{P,Q,T\}$. The power series $f$ verifies the Eggers condition at both $Q$ and $T$, but it does not satisfy this condition at the point $P$. Observe that $i_0(\Gamma_P,x)\neq i_0(f_P,x)-\de(P)$.
\end{enumerate}

\begin{figure}
    \begin{center}
\begin{tikzpicture}[scale=1]
   
    \draw [-, color=black,  line width=1pt](0,0) -- (0,1);
    \draw [->, color=black,  line width=1pt](0,1) -- (1,2.3);
    \draw [->, color=black,  line width=1pt](0,1) -- (-1,2.3);
    \draw [->, color=black,  line width=1pt](0,1) -- (0,0);
   
    \node [below, color=black] at (0,0) {$x$};
    \node [right, color=black] at (0,0.5) {\small{$1$}};
    \node [above, color=black] at (-1,2.3) {$Q$};
    \node [right, color=black] at (0.45,1.6) {\small{$p$}};
    \node [left, color=black] at (-0.45,1.6) {\small{$1$}};
    \node [above, color=black] at (1,2.3) {$T$};
                
    \node[draw,circle, inner sep=1.5pt,color=black, fill=black] at (0,1){};
    \node [left, color=black] at (0,1) {$P$}; 
    \end{tikzpicture}
\end{center}
 \caption{Eggers-Wall tree of Example \ref{ex:p&p+1}.}
\label{fig:ex}
   \end{figure}
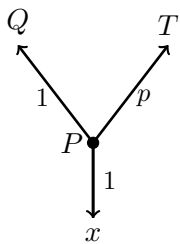
   \end{example}

\begin{example}\label{fig:exnewAP}
Let us come back to Example \ref{ex:pdoesnotimplyE}. Recall that $f$ does not satisfy the Eggers condition at the point $P$. Notice that the polar of $f$ is $\frac{\partial f}{\partial y }=y^{3}+x^{5}$ whose attaching point to $\Theta(f)$ is an unmarked point $P'\in(x,P)$ with $\ex(P')=\frac{5}{3}<\frac{5}{2}=\ex(P)$ (see the right-hand of Figure \ref{fig:ex2}). Observe that for any point $Q$ of the segment $(P',P]$ we have $i_0(\Gamma_Q,x)=0$ and $i_0(f_Q,x)-\de(Q)=2$.

\begin{figure}[h]
    \begin{center}
\begin{tikzpicture}[scale=1]
\begin{scope}[shift={(0,0)}]
   
    \draw [-, color=black,  line width=1pt](0,0) -- (0,1);
    \draw [->, color=black,  line width=1pt](0,1) -- (1,2.3);
    \draw [->, color=black,  line width=1pt](0,1) -- (-1,2.3);
    \draw [->, color=black,  line width=1pt](0,1) -- (0,0);
   
    \node [below, color=black] at (0,0) {$x$};
    \node [right, color=black] at (0,0.5) {\small{$1$}};
    \node [above, color=black] at (-1,2.3) {$Q$};
    \node [right, color=black] at (0.45,1.6) {\small{$2$}};
    \node [left, color=black] at (-0.45,1.6) {\small{$1$}};
    \node [above, color=black] at (1,2.3) {$T$};
                   
    \node[draw,circle, inner sep=1.5pt,color=black, fill=black] at (0,1){};
    \node [left, color=black] at (0,1) {$P$};
  \end{scope}

\begin{scope}[shift={(5,0)}]

    \draw [-, color=black,  line width=1pt](0,0) -- (0,1);
    \draw [->, color=black,  line width=1pt](0,1) -- (1,2.3);
    \draw [->, color=black,  line width=1pt](0,1) -- (-1,2.3);
    \draw [->, color=black,  line width=1pt](0,1) -- (0,0);
    \draw [->, color=blue,  line width=1pt](0,0.5) -- (-1,0.75);
   
    \node [below, color=black] at (0,0) {$x$};
    \node [right, color=black] at (0,0.25) {\small{$1$}};
    \node [right, color=black] at (0,0.75) {\small{$1$}};
    \node [left, color=blue] at (0,0.35) {\small{$P'$}};
    \node [above, color=black] at (-1,2.3) {$Q$};
    \node [right, color=black] at (0.45,1.6) {\small{$2$}};
    \node [left, color=black] at (-0.45,1.6) {\small{$1$}};
    \node [above, color=black] at (1,2.3) {$T$};
                                 
    \node[draw,circle, inner sep=1.5pt,color=black, fill=black] at (0,1){};
    \node[draw,circle, inner sep=1.5pt,color=blue, fill=blue] at(0,0.5){};
          
    \node [left, color=black] at (0,1) {$P$};
    \node [left, color=blue] at (-1,0.75) {$y^3+x^5$};

    \end{scope}          
    \end{tikzpicture}
\end{center}
 \caption{Eggers-Wall trees of Example \ref{fig:exnewAP}.}
\label{fig:ex2}
   \end{figure}
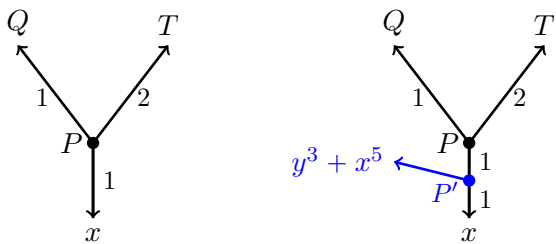
\end{example}

We are now in position to prove that the Eggers condition is sufficient for $\frac{\partial f}{\partial y}$ to admit Eggers decomposition.

\begin{corollary}\label{sufficient} [Eggers decomposition]\label{Eggersfact}
Let $f\in K[[x,y]]$ be such that $1<\ord f(0,y)<\infty$ and denote $\Gamma=\frac{\partial f}{\partial y}$. If $f$ satisfies Eggers condition, then $\frac{\partial f}{\partial y}$ admits Eggers decomposition. Moreover, if $f$ has no multiple irreducible factors, then for every leaf $P$ of $\Theta(f)$ we have $\Gamma^{(P)}=1$.
\end{corollary}

\begin{proof}
Assume that $f$ satisfies the Eggers condition. Then any pair $P',P$ of consecutive marked points of $\Theta(f)$ verifies the assumption of Theorem \ref{mainJ}. Hence the pair $(f,\Gamma)$ verifies property $(E2)$ for every point $Q\in\Theta(f)$ different from the root $R$. Furthermore, since $x$ is not a factor of $f$, we have $f_R=f_{P}=f$ where $P$ denotes the direct successor of $R$. Then by assumption $\ord f(0,y)=i_0(f,x)=i_0(f_P,x)\not\equiv 0 \pmod p$, so differentiating we get $i_0(\Gamma,x)=\ord\Gamma(0,y)=\ord f(0,y)-1=i_0(f,x)-1$. We conclude that $(f,\Gamma)$ also satisfies condition $(E2)$ for $Q=R$. Therefore $\frac{\partial f}{\partial y}$ admits Eggers decomposition (see Definition \ref{def:Ed}).

Finally, if all irreducible factors of $f$ are pairwise different and $P$ is a leaf of $\Theta(f)$, then $i_0(f_P,x)=\de(P)$. Since $(f,\Gamma)$ verifies condition $(E2)$, the pair $(f,\Gamma)$ also verifies condition $(E1)$. By item (ii) of $(E1)$, we get $i_0(\Gamma^{(P)},x)=0$ and we conclude that $\Gamma^{(P)}=1$.
\end{proof}

Notice that for $f$ irreducible the tree $\Theta(f)$ has only one leaf. The hypothesis $\ord f(0,y)\not\equiv 0 \pmod p$ in \cite[Theorem 1.5]{GB-P} is equivalent with the hypothesis of Corollary \ref{Eggersfact}, that is, $i_0(f_Q,x)\notequiv 0 \pmod p$ for any marked point $Q\in \Theta(f)$ different from the root. Indeed, for any marked point $Q\in \Theta(f)$ we have $f_Q=f$ and consequently $i_0(f_{Q},x) = i_0(f,x)=\ord f(0,y) \notequiv 0 \pmod p$, and $\frac{\partial f}{\partial y}$ admits Eggers decomposition. After Corollary \ref{coro:Edecg} and the description of the Eggers-Wall tree for $f$ irreducible (see section \ref{sect:E-Wtree}), we get $\frac{\partial f}{\partial y}=\prod_{i=1}^h \Gamma^{(P_i)}$ with 
 $i_0(\Gamma^{(P_i)},x)=n_1\cdots n_{i}-n_1\cdots n_{i-1}=n_1\cdots n_{i-1}(n_i-1)$ and $\frac{i_0(\Gamma^{(P_i)},f)}{i_0(\Gamma^{(P_i)},x)}=\frac{e_{i-1}\bb_i}{\bb_0}$. From Remark \ref{rem:computingd} we also get $\frac{i_0(g,f)}{i_0(g,x)}=\frac{e_{i-1}\bb_i}{\bb_0}$ for any irreducible factor $g$ of $\Gamma^{(P_i)}$. Hence we find the decomposition of  $\frac{\partial f}{\partial y}$ given in \cite[Theorem 1.5]{GB-P}, when $f$ is irreducible. Observe that even in the irreducible case, the Eggers condition does not imply a Puiseux factorization of the partial derivative, as the following example illustrates:

\begin{example}\label{ex:doesnotPuiseuxroot}
    Consider a power series $f(x,y)=y^{p+1}+x^3y^{p-1}+x^{p+1}y+x^p$ over a field of characteristic $p>2$. It is an irreducible Weierstrass polynomial of degree $p+1$, so $f$ satisfies the Eggers condition. Its partial derivative $\frac{\partial f}{\partial y}=y^p-x^3y^{p-2}+x^{p+1}$ is an irreducible Weierstrass polynomial which by \cite[Corollary 17.7]{Jaca} does not have any Puiseux root since it does not belong to $K[[x]][y^p]$.  
\end{example}

The following example illustrates the Eggers decomposition in the reducible case.

\begin{example}\label{ex:Eggersd}
Consider $f(x,y)=f_{1}f_{2}f_{3}\in K[[x,y]]$, where $K$ is an algebraically closed field of any characteristic different from $2$, $3$, and 11, $f_{1}=y^{2}+x^{3}$, $f_{2}=y^{3}+x^{4}$, and $f_{3}=(y^{3}+x^{4})^{2}+x^{9}$. 
Figure \ref{fig:example} shows the Eggers-Wall tree of $f$. We get
\begin{itemize}
\item $\ex(P)=\frac{4}{3}$, $\ex(Q)=\frac{3}{2}$, and $\ex(T)=\frac{11}{6}$,
\item $\ic(P)=\ex(P)$, $\ic(Q)=\ic(P)+\frac{1}{1}(\ex(Q)-\ex(P))=\ex(Q)$ and $\ic(T)=\ic(P)+\frac{1}{3}(\ex(T)-\ex(P))=\frac{3}{2}$.
\end{itemize}

Notice that $f_{P}=f$, $f_{Q}=f_{1}$, and $f_{T}=f_{2}f_{3}$. Hence $\Gamma=\frac{\partial f}{\partial y}=\Gamma^{(P)}\Gamma^{(Q)}\Gamma^{(T)}$, where
$i_{0}(\Gamma^{(P)},x)=\de(Q)+\de(T)-\de(P)=1+3-1=3$, $i_{0}(\Gamma^{(Q)},x)=2-1=1$, and $i_{0}(\Gamma^{(T)},x)=6+3-3=6$. Moreover, by Remark \ref{rem:computingd} we can compute $d(f_i,h)$ for any irreducible factor $h$ of $\Gamma$ and $f_i$ for $i\in\{1,2,3\}$. For example if $h$ is an irreducible factor of $\Gamma^{(T)}$, then
\[
d(f_i,h) =\left\{\begin{array}{ll} 
\ic(P)&\;\hbox{\rm if $i=1$}\\
\ic(T) &\;\hbox{\rm if $i\in \{2,3\}$.}
\end{array}
\right.
\]

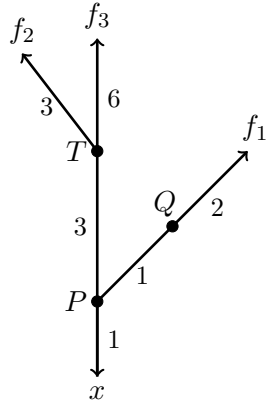
\begin{figure}
    \begin{center}
\begin{tikzpicture}[scale=1]
   
    \draw [-, color=black,  line width=1pt](0,0) -- (0,1);
    \draw[->, color=black,  line width=1pt](0,1) -- (0,4.5);
    \draw [->, color=black,  line width=1pt](0,1) -- (2,3);
    \draw [->, color=black,  line width=1pt](0,3) -- (-1,4.3);
    \draw [->, color=black,  line width=1pt](0,1) -- (0,0);
   
    \node [below, color=black] at (0,0) {$x$};
    \node [right, color=black] at (0,0.5) {\small{$1$}};
    \node [above, color=black] at (-1,4.3) {$f_{2}$};
    \node [below, color=black] at (0.6,1.6) {\small{$1$}};
    \node [below, color=black] at (1.6,2.5) {\small{$2$}};
    \node [left, color=black] at (-0.45,3.6) {\small{$3$}};
    \node [above, color=black] at (2.1,3) {$f_{1}$};
    \node [above, color=black] at (0,4.5) {$f_{3}$};
    \node [right, color=black] at (0,3.7) {\small{$6$}};
    \node [left, color=black] at (0,2) {\small{$3$}};
                   
    \node[draw,circle, inner sep=1.5pt,color=black, fill=black] at (0,1){};
    \node[draw,circle, inner sep=1.5pt,color=black, fill=black] at (0,3){};
    \node[draw,circle, inner sep=1.5pt,color=black, fill=black] at (1,2){};
          
    \node [left, color=black] at (0,1) {$P$}; 
    \node [left, color=black] at (0,3) {$T$};
    \node [above, color=black] at (0.9,2) {$Q$};  
    \end{tikzpicture}
\end{center}
 \caption{Eggers-Wall tree of $f=f_1f_2f_3$ of Example \ref{ex:Eggersd}.} 
\label{fig:example}
   \end{figure}

\end{example}

\section{Eggers condition is necessary for Eggers decomposition}\label{sec:5}

Let $f\in K[[x,y]]$ be such that $1<\ord f(0,y)<\infty$ and assume that $\frac{\partial f}{\partial y}$ admits Eggers decomposition. The goal of this section is to show that $f$ satisfies the Eggers condition. We begin by proving that $f$ satisfies $\de$--condition. Our strategy relies on the substitution of $x$ by $x^{p^k}$ in the power series $f$, which preserves the shape of the Eggers-Wall tree (i.e., how branches split remains unchanged after the substitution). This section starts with a few lemmas required to establish this invariance.

\begin{lemma}\label{L:p}
For every irreducible Weierstrass polynomial $f\in K[[x]][y]$ and positive integer $k$ there exist an irreducible Weierstrass polynomial $\bar f\in K[[x]][y]$ and a divisor $s$ of $p^k$ such that $f(x^{p^k},y)=\bar f(x,y)^s$.
\end{lemma}

\begin{proof}
    Let $f=\prod_{j=1}^s  f_j$ be the factorization of $f$ into irreducible monic polynomials in $K[[x^{1/p^k}]][y]$. The characteristic of $K$ equals $p$, so we have that ${ f_j}^{p^k}\in K[[x]][y]$ for $1\leq j\leq s$.  Since $K[[x]][y]$ is a unique factorization domain and $f$ is irreducible, from the equality $f^{p^k}=\prod_{j=1}^s { f_j}^{p^k}$ we deduce that for every $j\in\set{1,\ldots,s}$ there exists a positive integer $r_j$ such that ${ f_j}^{p^k}=f^{r_j}$. Observe that $p^k=\sum_{j=1}^s{r_j}$. Then, for $i\neq j$, we get $ f_j^{r_i p^k}=f^{r_j r_i}= f_i^{r_j p^k}$. Since $K[[x^{1/p^k}]][y]$ is a unique factorization domain and $f_i$ and $ f_j$ are irreducible, we must have $f_i= f_j$ and $r_i=r_j$. Hence, $f={ f_1}^s$ with $p^k=sr_1$. To end the proof, it is enough to take $\bar f(x,y)= f_1(x^{p^k},y)$.
\end{proof}

\begin{lemma}\label{L:p1}
Let $\alpha(x)$ be a Puiseux series of positive order and let $k$ be a positive integer. Let the Weierstrass polynomials $f,\hat f  \in K[[x]][y]$ be the minimal polynomials of $\alpha(x)$, $\alpha(x^{p^k})$ respectively. Then $f(x^{p^k},y)= \hat f(x,y)^s$ for some positive integer $s$. In addition, if $p^k$ divides the index of $\alpha(x)$ then $s=p^k$. If $p$ does not divide the index of $\alpha(x)$, then $s=1$.
\end{lemma}

\begin{proof}
By Lemma \ref{L:p}, there exist an irreducible Weierstrass polynomial $\bar f\in K[[x]][y]$ 
and a divisor $s$ of $p^k$ such that $f(x^{p^k},y)=\bar f(x,y)^s$. Replacing $y$ with $\alpha(x^{p^k})$ and taking into account that $\alpha(x)\in\Zer(f)$, we get $\bar f(x,\alpha(x^{p^k}))=0$. Hence $\hat f=\bar f$. This proves the first part of the lemma. 

Observe that $\deg f=s\deg\bar f$. Moreover, by Property \ref{prop:minpol} the index of $\alpha(x)$ (resp. of $\alpha(x^{p^k})$) equals $\deg f$ (resp. $\deg\bar f$). If $dp^k$ is the index of $\alpha(x)$, then the index of $\alpha(x^{p^k})$ is equal to $d$, and this leads to $s=p^k$. If $p$ does not divide the index of $\alpha(x)$, then the index of $\alpha(x^{p^k})$ is equal to that of $\alpha(x)$, and thus $s=1$.
\end{proof}

\begin{lemma}\label{lem:dis}
Let $f,\bar f, g, \bar g\in K[[x]][y]$ be irreducible Weierstrass polynomials such that 
$f(x^n,y)=\bar f(x,y)^a$ and $g(x^n,y)=\bar g(x,y)^b$ for some positive integers $a,b,n$. Then
$d(\bar f, \bar g)=n d(f,g)$.
\end{lemma}
\begin{proof} We have
\begin{eqnarray*}
n\,d(f,g)&=& \frac{n\,i_0(f(x,y),g(x,y))}{i_0(f(x,y),x)\,i_0(g(x,y),x)}=\frac{i_0(f(x^n,y),g(x^n,y))}{i_0(f(x^n,y),x)\,i_0(g(x^n,y),x)}\\
&=&\frac{i_0(\bar f^a,\bar g^b)}{i_0(\bar f^a,x)\,i_0(\bar g^b,x)}=\frac{ab\,i_0(\bar f,\bar g)}{ai_0(\bar f,x)\,bi_0(\bar g,x)}=d(\bar f,\bar g),
\end{eqnarray*}

\noindent where we use the well-known formula 
$i_0(f(x^n,y), g(x^n,y))=ni_0(f(x,y),g(x,y))$ (see \cite[Proposition 3.12, Chapter 2]{P-P} 
\end{proof}

\begin{remark}\label{rem:hat}
Write a power series $f\in K[[x,y]]$ such that $1< \ord f(0,y)<\infty$ as the 
product $u\prod_{i=1}^r f_i^{a_i}$ where $u$ is a unit in $K[[x,y]]$ and 
$f_i\in K[[x]][y]$ are pairwise coprime irreducible Weierstrass polynomials. Let $k$ be a positive integer. By Lemma~\ref{L:p} there exist irreducible Weierstrass polynomials $\bar f_i\in K[[x]][y]$ and positive integers $s_i$ 
such that $f_i(x^{p^k},y)=\bar f_i(x,y)^{s_i}$ for $i=1,\dots, r$.  Consider any power series $\check f$ of the form $\bar u\prod_{i=1}^r \bar f_i^{b_i}$, where $\bar u$ is a unit in $K[[x,y]]$ and $b_i$ are positive integers for $i=1,\dots, r$. In particular $f(x^{p^k},y)$ has this form. Note that $f_i\to \bar f_i$ defines a one-to-one correspondence between the leaves of $\Theta(f)$ and the leaves of $\Theta(\check f)$. In addition, by Lemma~\ref{lem:dis} we obtain 
$d(\bar f_i, \bar f_j)=p^k d(f_i,f_j)$ for $1\leq i<j\leq r$. The inequalities between the logarithmic distances $d(f_i,f_j)$ determine the distribution of 
ramification points in $\Theta(f)$. Thus the shapes of $\Theta(f)$ and $\Theta(\check f)$ (where we neglect bamboo points) are the same. Moreover for any ramification point $P=\langle x,f_i,f_j\rangle$ of $\Theta(f)$ 
and its corresponding ramification point $\bar P=\langle x,\bar f_i,\bar f_j\rangle$
of $\Theta(\check f)$, we have $\ic(\bar P)=p^k \ic(P)$.
\end{remark}

The following result will be key in showing that the Eggers condition is necessary for Eggers decomposition.

\begin{pro}\label{i-condition}
Let $f\in K[[x,y]]$ be such that $1<\ord f(0,y)<\infty$. If $\frac{\partial f}{\partial y}$ admits Eggers decomposition then $f$ satisfies $\de$--condition.
\end{pro}

\begin{proof}
Assume that the conclusion of the proposition does not hold. Then there exist consecutive marked 
points $P$, $P'$ of $\Theta(f)$ such that $\de(P) \notequiv 0 \pmod p$ and $\de(P') \equiv 0 \pmod p$.  
Let $f_j$ be an irreducible factor of $f_{P'}$. Since $P$ is a point of discontinuity of the restriction of the index function to $[x,f_j]$, for some 
key polynomial $g\in K[[x]][y]$ of $f_j$ we have $d(g,f_j)=\ic(P)$. The degree of $g$ equals $\de(P)$, which is not divisible by $p$, so there exists a Puiseux root $\alpha(x)$ of $g$ of index $\de(P)$. Observe that $g$ is the minimal polynomial of $\alpha(x)$.  
As $f_j$, the Puiseux series $\alpha(x)$, and $g$ satisfy the assumptions of Corollary \ref{C2}, 
there exists a Puiseux series $\beta(x)=\alpha(x)+cx^{\ex(P)}$ such that $d(f_j,g_\beta)>\ic(P)$, where $g_\beta$ denotes the minimal polynomial of $\beta(x)$. 
Then we have $P\prec_x P_{g_\beta}$ in $\Theta(f)$, and therefore the degree of $g_\beta$ is a multiple of $\de(P')$, and also of $p$.

Denote $h:=\ex(P)$. We can write $h=\frac{m}{q\deg g}$ with $\gcd(m,q)=1$. Indeed, if $h=a/b$ with $\gcd(a,b)=1$ and $D=\gcd(b,\deg g)$, then $m=a\deg g/D$ and $q=b/D$.
Since we have $\deg g_\beta=\hbox{index of}\,\beta(x)=\hbox{lcm}(\deg g,b)=q\deg g$, $\deg g_\beta\equiv 0 \pmod p$, and $\deg g\notequiv 0 \pmod p$, we deduce that $q$ is divisible by $p$. 
Hence we can write $h=\frac{m}{dp^k \deg g}$, where $\gcd(m,dp^k)=\gcd(d,p)=1$ and $k\geq1$. Moreover, $\deg g_\beta=dp^k \deg g$.

Let $b_0=\deg g$, $b_1,\dots, b_{\ell}$ be the Puiseux characteristic sequence of $g$. For $i=1,\dots,\ell$, denote $h_i=b_i/b_0$. Note that by construction, $h$ is bigger than all the $h_i$. 

\medskip
{\it Claim 1:} $g$ is the last by one key polynomial of $g_\beta$.

Indeed, for $i \in \{1,\dots,\ell\}$, let $g^{(i-1)}$ be the minimal polynomial 
of the truncation of $\alpha(x)$ consisting of the sum of the terms of exponents strictly less than $h_i$. Then by \cite[Proposition 4.4]{P-P}, $g^{(0)},\ldots, g^{(\ell-1)}, g^{(\ell)}=g$ is a sequence of key polynomials of $g$. By \cite[Theorem 5.2]{GB-P-2015}, since $d(f_j,g_\beta)>\ic(P)$ and $g$ is a key polynomial of $f_j$, the sequence $\{ g^{(i)}\}_{0\leq i\leq \ell}$
is the initial part of a sequence of key polynomials of $g_\beta$. Set $\bar b_i=i_0(g_\beta,g^{(i-1)})$ for $i=1,\dots,\ell+1$. To prove the claim, it is enough to show that 
\[ \gcd(\deg g_\beta,\bar b_1,\ldots \bar b_{\ell+1})=1.\]

Let $\bar g$ be the minimal polynomial of $\alpha(x^{p^k})$. Since the index of $\alpha(x)$ is not divisible by $p$, the sequence $b_0,p^kb_1,\dots, p^kb_{\ell}$ is the Puiseux characteristic sequence of $\bar g$. 
Thus by Lemma~\ref{L:p1}, $\bar g^{(i-1)}(x,y) = g^{(i-1)}(x^{p^k},y)$, $1\leq i\leq\ell$, are the minimal polynomials of the corresponding truncations of $\alpha(x^{p^k})$, and $\bar g(x,y)=g(x^{p^k},y)$. It follows that $\bar g^{(0)},\ldots,\bar g^{(\ell)}=\bar g$ is a sequence of key polynomials of $\bar g$. 

Let $\bar g_{\beta}$ be the minimal polynomial of $\beta(x^{p^k})$. The index of $\beta(x)$ is divisible by $p^k$, hence by Lemma \ref{L:p1}, we have $g_\beta(x^{p^k},y)=\bar g_\beta(x,y)^{p^k}$. This implies that $\deg \bar g_\beta=\deg g_\beta/p^k=d\deg g$.  Note that $k(\bar g_\beta,\bar g)=p^kh$ and thus by Property 2.3, the exponent function takes the value $p^kh$ at the  point of the tree $\Theta(\bar g_\beta \bar g)$. Next we compute, from its definition, the integers following $d\deg g$ in the Puiseux characteristic sequence of $\bar g_\beta$, and finally we compare the truncations of $\alpha(x^{p^k})$ and $\beta(x^{p^k})$ to obtain  key polynomials of $\bar g_\beta$. Recall that $h=\frac{m}{dp^k\deg g}$ where $\gcd(m,p^kd)=1$. We obtain the following:
\begin{itemize}
    \item If $d=1$, then $\deg g$, $p^kb_1,\dots, p^kb_{\ell}$ is the Puiseux characteristic sequence of $\bar g_\beta$ and $\bar g^{(0)},\ldots,\bar g^{(\ell-1)}$, $\bar g_\beta$ is a sequence of key polynomials of $\bar g_\beta$.
    \item If $d>1$, then $d \deg g$, $p^kdb_1,\dots, p^kdb_{\ell},m$ is the Puiseux characteristic sequence of $\bar g_\beta$. Hence  
$\bar g^{(0)},\ldots,\bar g^{(\ell-1)}$ is the initial part of a sequence of key polynomials of $\bar g_\beta$. Now let $\bar g_1$ be the minimal polynomial of the truncation of $\beta(x^{p^k})$ consisting of the terms of exponents strictly less than $p^kh$. Then $\bar g_1$ is an $\ell$-key polynomial of $\bar g_\beta$. Since $\deg \bar g_1=\deg \bar g$ and $k(\bar g_\beta,\bar g_1)=k(\bar g_\beta,\bar g)$, we get that $\bar g$ is also an $\ell$-key polynomial of $\bar g_\beta$. Hence $\bar g^{(0)},\ldots,\bar g^{(\ell-1)},\bar g^{(\ell)}=\bar g,\bar g_\beta$ is a sequence of key polynomials of $\bar g_\beta$.
\end{itemize} 

\noindent In both cases 
\begin{equation}\label{semigroup}
\gcd(\, \deg \bar g_\beta,i_0(\bar g_{\beta},\bar g^{(0)}),\dots, i_0(\bar g_{\beta},\bar g^{(\ell)}) \,) =1.
\end{equation}

By equalities $d(\bar g_\beta,\bar g^{(i)})=p^kd(g_\beta,g^{(i)})$ (see Lemma \ref{lem:dis}), $\deg g_\beta=p^k \deg\bar g_\beta$, and $\deg g^{(i)}=\deg \bar g^{(i)}$, we get 
\[ \bar b_i=i_0(\bar g_\beta,\bar g^{(i-1)})\quad \mbox{for } i=1,\dots,\ell+1. \]
Since $d(g_{\beta},g^{(\ell-1)})= d(g,g^{(\ell-1)})$, we have that 
$i_0(g_\beta,g^{(\ell-1)})=i_0(g,g^{(\ell-1)}) (\deg g_\beta / \deg g)$, thus $\bar b_\ell \equiv 0 \pmod p$. By the integral formula for the exponent function of the Eggers-Wall tree $\Theta(\bar g_\beta \bar g)$, we get $\deg \bar g\bigl(d(\bar g_{\beta},\bar g)-d(\bar g_{\beta},\bar g^{(\ell-1)})\bigr)=p^kh-p^kh_{\ell}$.
Rewriting this formula in terms of intersection multiplicities we obtain 
$\bar b_{\ell+1}-(\deg \bar g/ \deg \bar g^{(\ell-1)} )\bar b_{\ell}=\deg \bar g_{\beta}(p^kh-p^kh_{\ell})=m-p^k h_\ell \deg \bar g_\beta$. 
After reduction modulo $p$,
\[ \bar b_{\ell+1} \equiv m \pmod p.\]
Finally by~(\ref{semigroup}) we get 
$\gcd(\deg g_\beta, \bar b_1,\ldots \bar b_{\ell+1})=
\gcd(p^k\deg \bar g_\beta,\bar b_1,\ldots \bar b_{\ell+1})=1$
because $p$ does not divide $\bar b_{\ell+1}$. Claim 1 is proved. 

\medskip

Let $\gamma(x)=\alpha(x)+cx^{h}+x^{h'}$, where $h'=\frac{m'}{d'}>h$ and $d'$ is coprime with the index of $\beta(x)$. 
Let $g_{\gamma}$ be the minimal polynomial of $\gamma(x)$. Note that $P$ can also be seen as a point of $\Theta(g_\beta)$ and it is the attaching point of $g$ to this tree.

\medskip
\textit{Claim 2:} If $Q$ is the attaching point of $g_{\gamma}$ to $\Theta(g_{\beta})$, then 
\begin{equation}\label{claim2}
e(Q)-e(P)=p^k(h'-h) .
\end{equation}

Let $\bar g_{\gamma}$ be the minimal polynomial of $\gamma(x^{p^k})$. Since the index of $\gamma(x)$ is divisible by $p^k$, it follows from Lemma \ref{L:p1} that 
$g_\gamma(x^{p^k},y)=\bar g_\gamma(x,y)^{p^k}$. In addition, by construction the index of $\gamma(x^{p^k})$ is not divisible by $p$, thus $\gcd(\deg\bar g_\gamma,p)=1$.

Let $\bar P$ and $\bar Q$ be the attaching points of $\bar g$ and $\bar g_\gamma$, respectively, to the Eggers-Wall tree $\Theta(\bar g_\beta)$. By abuse of notation, we use the bold letters $\ic$, $\de$, and $\ex$ to denote the functions defined on both $\Theta(g_\beta)$ and $\Theta(\bar g_\beta)$. We know already that $\ex(\bar P)=p^k h$ (observe that $\bar P$ can be seen as the  point of the tree $\Theta(\bar g_\beta \bar g)$). Property \ref{prop:1} also allows us to compute $\ex(\bar Q)$. Indeed, $\ex(\bar Q)=k(\bar g_\gamma,\bar g_\beta)=\ord(\gamma(x^{p^k})-\beta(x^{p^k}))=p^k h'$. We get 
\begin{equation}\label{eq:difexp}
    p^k(h'-h)=\ex(\bar Q)-\ex(\bar P),
\end{equation}
and therefore $\ex(\bar P)<\ex(\bar Q)$ and also $\ic(\bar P)<\ic(\bar Q)$. This last inequality and the fact that 
\begin{equation}\label{eq:contactpk}
\ic(\bar P)=p^k\ic(P)\hbox{ and }\ic(\bar Q)=p^k\ic(Q)
\end{equation}
(see Lemma \ref{lem:dis}) imply that $\ic(P)<\ic(Q)$.

By Claim 1, $g$ is the last by one key polynomial of $g_\beta$, hence for all $A\in(P,g_\beta]$, we have $\de(A)=\deg{g_\beta}$. Moreover, we know also that $\bar B\preceq_x\bar P$ in $\Theta(\bar g_\beta)$, where $\bar B$ is the attaching point of the last by one key polynomial of $\bar g_\beta$ to $\Theta(\bar g_\beta)$. As a consequence, $\de(\bar A)=\deg{\bar g_\beta}$ for all $\bar A\in(\bar P,\bar g_\beta]$. Then, 
\begin{align*}
\ex(\bar Q)-\ex(\bar P)&=
\deg\bar g_\beta[\ic(\bar Q)-\ic(\bar P)] = 
p^k\deg \bar g_\beta[\ic(Q)-\ic(P)] \\
&=\deg g_\beta[\ic(Q)-\ic(P)] = \ex(Q)-\ex(P),
\end{align*}
where the first and the fourth equalities follow from the integral formula for $\ex$, the second equality follows from \eqref{eq:contactpk}, and the third equality follows from the relation $\deg g_\beta=p^k\deg \bar g_\beta$. Taking into account equation \eqref{eq:difexp}, this ends the proof of Claim 2.

\medskip

Let $ f^{*}(x,y)=f(x^{p^k},y)$. Applying Remark \ref{rem:hat} to the power series $f g g_\beta$ and $ f^{*} \bar g \bar g_\beta$, we have a one-to-one correspondence between the ramification points of $\Theta(f g g_\beta)$ and those of $\Theta(f^{*} \bar g \bar g_\beta)$. The points $P$ and $P_{g_\beta}$ can be seen in the tree $\Theta(f g g_\beta)$ and they are both ramification points. Note that the ramification point corresponding to $P$ coincides with the already defined $\bar P$, seen in $\Theta(f^{*} \bar g \bar g_\beta)$, and that $\bar P\prec_x\bar P_{g_\beta}$, where $\bar P_{g_\beta}$ is the ramification point corresponding to $P_{g_\beta}$. Let us distinguish two cases:
\begin{itemize}
    \item $P'\preceq_x P_{g_\beta}$. Since $P'$ is a marked point and the index function is constant in $(P,g_\beta]$, in this case the point $P'$ is necessarily a ramification point. There are no marked points in $(P,P')$ and the index function is constant in $(\bar P,\bar g_\beta]$, therefore we conclude that there are no marked points in $(\bar P,\bar P')$, where $\bar P'$ is the ramification point corresponding to $P'$. Set $T=P'$ and $\bar T=\bar P'$. 
    \item $P_{g_\beta}\prec_x P'$. The same arguments as before allow us to conclude that there are no marked points in $(\bar P,\bar P_{g_\beta})$. Set $T=P_{g_\beta}$ and $\bar T=\bar P_{g_\beta}$. 
\end{itemize}

Take $h'$ such that  $h'-h$ is sufficiently close to $0$. Then by~\eqref{claim2} the attaching point $Q$ of $g_\gamma$ to the Eggers-Wall tree of $f$ belongs to $(P,T)$ and by~\eqref{eq:difexp} the attaching point $\bar Q$ of $\bar g_\gamma$ to the Eggers-Wall tree of $ f^{*}$ belongs to $(\bar P,\bar T)$. Note that $\de(\bar T)=\deg \bar g_\beta\notequiv 0 \pmod p$. Hence the segment $[\bar P,\bar T]$ (seen in the tree $\Theta(f^{*})$), the rational number $\iota=p^k h'$, and the Puiseux series $\gamma(x^{p^k})$ satisfy the assumptions of Lemma~\ref{edgeE23}.

It follows that the Newton polygon of $\widetilde{f^{*}}(x,y)= f^{*}(x,y+\gamma(x^{p^k}))$ has an edge $E$ of inclination $\iota$ that touches the horizontal axis. Moreover, the Newton principal part of $\widetilde{f^{*}}$ with respect to $E$ is of the form $\epsilon x^a(y-x^\iota)^\tau$, with $\tau=i_0(f^{*}_{\bar T},x)/\de(\bar T)$, $\epsilon\in K^*$, and $a\geq0$. We have $\ord \widetilde {f^{*}}(x,0)=a+\iota\tau$. Since the point $(a+\iota(\tau-1),1)$ belongs to $E$, we have that $\ord \frac{\partial \widetilde {f^{*}}}{\partial y}(x,0)\geq a+\iota(\tau-1)$. Therefore we get 
\[\ord f^{*}(x,\gamma(x^{p^k}))-\ord \frac{\partial f^{*}}{\partial y}(x,\gamma(x^{p^k})) = 
\ord \widetilde {f^{*}}(x,0)-\ord \frac{\partial \widetilde {f^{*}}}{\partial y}(x,0) \leq p^kh',\]
and consequently $\ord  f(x,\gamma(x))-\ord \frac{\partial f}{\partial y}(x,\gamma(x)) \leq h' = h+(h'-h) < \ex(P)+p^k(h'-h) = \ex(Q)$, which contradicts Corollary~\ref{co:pcondition}.
\end{proof}

Note that the statement of Proposition \ref{i-condition} is not an equivalence. Indeed, if $f$ is the power series of Example \ref{ex:pdoesnotimplyE}, then $\de(A)$ is not a multiple of $3$ for any marked point $A\in\Theta(f)$. However, $\frac{\partial f}{\partial y}$ does not admit Eggers decomposition (see Example \ref{fig:exnewAP}).

We are now ready to prove the main result of the section.

\begin{theorem}\label{necessary}
    Let $f\in K[[x,y]]$ be such that $1<\ord f(0,y)<\infty$. If $\frac{\partial f}{\partial y}$ admits Eggers decomposition, then $f$ satisfies Eggers condition.
\end{theorem}

\begin{proof}
Assume that the pair $(f,\frac{\partial f}{\partial y})$ satisfies condition (E2). Let $P$ be a marked point of $\Theta(f)$ different from the root
and let $[P',P]$ be a segment of $\Theta(f)$ with no marked points in $(P',P)$. By Proposition~\ref{i-condition}, $\de(P)\notequiv 0 \pmod p$.
Let $\iota=m/n$ be a rational number such that ${\rm gcd}(m,n)={\rm gcd}(p\cdot\de(P),n)=1$ and $\ex(P')<\iota<\ex(P)$. 
Then there exists a Puiseux series $\alpha(x)$ satisfying the conclusions of Lemma~\ref{edgeE23}. In particular, 
the Newton diagram of $\tilde f(x,y)=f(x,y+\alpha(x))$ has a compact edge $E$ such that 
$\tilde f_{\vert E}=\varepsilon x^a(y-x^\iota)^{i_0(f_P,x)/\de(P)}$, with $\varepsilon\in K^*$ and $a\geq0$. On the one hand, by Corollary~\ref{co:pcondition}, we have  $\ord \tilde f(x,0)-\ord \frac{\partial \tilde f}{\partial y}(x,0)=\ord f(x,\alpha(x))-\ord \frac{\partial f}{\partial y}(x,\alpha(x))=\iota.$ On the other hand, by the form of $\tilde f_{\vert E}$ we get $\ord\tilde f(x,0)-\ord \frac{\partial \tilde f}{\partial y}(x,0) = \iota$  
if and only if $\frac{\partial \tilde f_{\vert E}}{\partial y}(x,0)\neq 0$. It follows that 
$\frac{\partial \tilde f_{\vert E}}{\partial y} \neq 0$, which implies that $i_0(f_P,x)/\de(P) \notequiv 0 \pmod p$. 
Hence $i_0(f_P,x)\notequiv 0 \pmod p$, that is, $f$ satisfies the Eggers condition at the point $P$.
\end{proof}

After Theorem \ref{necessary} and Corollary \ref{sufficient} we conclude that Eggers condition is a necessary and sufficient condition for Eggers decomposition of the partial derivative. 

\bigskip

\noindent{\bf Funding information:} The first-named author was partially supported by the AGAUR project 2021 SGR 00603.
The second-named author was partially supported by the Spanish grant PID2023-149508NB-I00 funded by MICIU/AEI/10.13039/501100011033 and by FEDER, UE.
The three authors were also supported by the Spanish grant PID2019-105896-GB-I00 funded by MCIN/AEI/10.13039/501100011033. Part of this work was carried out during a visit of the third-named author to Universitat Polit\`ecnica de Catalunya - BarcelonaTech (UPC), partially funded by project ALECTORS 2023 - R02362.

\vspace{0.5cm}

\noindent {Ana Bel\'en de Felipe}\\
	Department of Mathematics, Barcelona East School of Engineering (EEBE), Universitat Polit\`ecnica de Catalunya - BarcelonaTech (UPC), Campus Diagonal Bes\`os, Ed. A, Avda. Eduard Maristany, 16, 08019 Barcelona, Spain. \\ORCID ID: 0000-0002-2171-9829 
	
	 \noindent {ana.belen.de.felipe@upc.edu}
	\vspace{0.3cm}

\noindent {Evelia Rosa Garc\'{\i}a Barroso}\\
	Departamento de Matem\'aticas, Estad\'{\i}stica e I.O.\\
	Instituto Universitario de Matem\'aticas y Aplicaciones (IMAULL). Universidad de La Laguna.\\
	Apartado de Correos 456.
	38200 La Laguna, Tenerife, Espa\~na \\ORCID ID: 0000-0001-7575-2619
	
	 \noindent {ergarcia@ull.es}
	\vspace{0.3cm}

\noindent {Janusz Gwo\'zdziewicz}\\
	Institute of Mathematics, University of the National Education Commission, Krakow \\
    Podchor\c{a}\.{z}ych 2 PL-30-084 Cracow, Poland \\ORCID ID: 0000-0002-2035-1073
	
	 \noindent {janusz.gwozdziewicz@uken.krakow.pl}

\end{document}